\documentclass[a4paper]{amsart}
\usepackage{a4wide}
\usepackage[all]{xy}
\usepackage[colorlinks=true]{hyperref}
\usepackage{enumerate}
\usepackage{amssymb}

\newtheorem{thm}{Theorem}[section]
\newtheorem{prop}[thm]{Proposition}
\newtheorem{lem}[thm]{Lemma}
\newtheorem{cor}[thm]{Corollary}

\theoremstyle{definition}
\newtheorem{definition}[thm]{Definition}
\newtheorem{example}[thm]{Example}

\theoremstyle{remark}
\newtheorem{remark}[thm]{Remark}
\newtheorem{remarks}[thm]{Remarks}

\newenvironment{ex}[1]{\begin{example}#1}{\[ \Diamond \]\end{example}}
\newenvironment{defn}[1]{\begin{definition}#1}{\[ \Diamond \]\end{definition}}
\newenvironment{rem}[1]{\begin{remark}#1}{\[ \Diamond \]\end{remark}}
\newenvironment{rems}[1]{\begin{remarks}#1}{\begin{center}$\Diamond$ \end{center}\end{remarks}}
\newcommand{\Rr}{\mathbb R}
\newcommand{\Ss}{\mathbb S}
\newcommand{\Zz}{\mathbb Z}

\newcommand{\Cc}{\mathbb C}
\newcommand{\CP}{\mathbb CP}
\newcommand{\Tt}{\mathbb T}
\renewcommand{\d}{\mathrm d}

\newcommand{\set}[1]{\left\{#1\right\}}

\newcommand{\eps}{\varepsilon}

\newcommand{\X}{\ensuremath{\mathfrak{X}}}
\newcommand{\F}{\ensuremath{\mathcal{F}}}

\newcommand{\Lie}{\boldsymbol{\pounds}}    % Lie derivative
\newcommand{\T}{\mathcal T}                % Generalized tangent bundle
\DeclareMathOperator{\Ad}{Ad}           % Adjoint
\DeclareMathOperator{\red}{red}         % Reduced space
\DeclareMathOperator{\Ker}{Ker}         % Kernel
\newcommand{\al}{\alpha}                
\newcommand{\be}{\beta}

\newcommand{\G}{\mathcal{G}}            % Lie groupoid
\newcommand{\s}{\mathbf{s}}             % source map
\renewcommand{\t}{\mathbf{t}}           % target map
\renewcommand{\gg}{\mathfrak{g}}        % Lie algebra
\newcommand{\hh}{\mathfrak{h}}          % Lie subalgebra
\newcommand{\pp}{\mathfrak{p}}        % Factor in Cartan decomposition

\newcommand{\tto}{\rightrightarrows}    % Arrows of a groupoid

\newcommand{\leaf}{\mathcal S}          % Symplectic leaf

\usepackage[dvips]{graphicx}

\begin{document}

\title{The momentum map in Poisson geometry}
\author{Rui Loja Fernandes}
\address{Departamento de Matem\'{a}tica\\
Instituto Superior T\'{e}cnico\\1049-001 Lisboa\\ Portugal} 
\email{rfern@math.ist.utl.pt}
\author{Juan-Pablo Ortega}
\address{Centre National de la Recherche Scientifique\\
Universit\'e de Franche-Comt\'e\\ F-25030 Besan\c{c}on cedex\\ France}
\email{Juan-Pablo.Ortega@univ-fcomte.fr}
\author{Tudor S. Ratiu}
\address{Section de Math\'ematiques and Bernoulli Center\\ \'Ecole Polytechnique F\'ed\'erale
de Lausanne\\ CH-1015 Lausanne\\ Switzerland.}
\email{Tudor.Ratiu@epfl.ch}
\thanks{RLF was supported by the Funda\c{c}\~ao para a Ci\^encia e a Tecnologia
through the Program POCI 2010/FEDER and by the Projects POCI/MAT/57888/2004 and POCI/MAT/55958/2004. JPO was partially supported by a "Bonus Qualit\'e Recherche" of the Universit\'e de Franche-Comt\'e and by the Agence National de
la Recherche, contract number JC05-41465. TSR was partially supported by the Swiss National Science Foundation grant 200020-109054.}

%\date{April 28, 2007}

\begin{abstract}
Every action on a Poisson manifold by Poisson diffeomorphisms lifts to a 
Hamiltonian action on its symplectic groupoid which has a canonically 
defined momentum map. We study various properties of this momentum map as well as its use in reduction.
\end{abstract}

%%% ----------------------------------------------------------------------
\maketitle
%%% ----------------------------------------------------------------------

%%%%%%%%%%%%%%%%%%%%%%%%%%%%%%%%%%%%
%%%%%%%%%%%%%%%%%%%%%%%%%%%%%%%%%%%%
%%%%%%%%%%%%%%%%%%%%%%%%%%%%%%%%%%%%

\section{Introduction}             %
\label{sec:introduction}           %
%%%%%%%%%%%%%%%%%%%%%%%%%%%%%%%%%%%%
%%%%%%%%%%%%%%%%%%%%%%%%%%%%%%%%%%%%
%%%%%%%%%%%%%%%%%%%%%%%%%%%%%%%%%%%%

Let $\Phi:G\times M\to M$ be a smooth action of a Lie group $G$ on a
Poisson manifold $(M,\{\cdot,\cdot\})$. Let $\mathfrak{g} $ be the Lie algebra of $G$ and $\mathfrak{g}^\ast$ its dual. We assume that this is a \emph{Poisson action}
i.e., $G$ acts by Poisson diffeomorphisms. Such a Poisson action 
will not usually have a momentum map in the classical sense. For example, 
there can be no momentum map if the action maps points in a symplectic leaf to
points in a different symplectic leaf. However, one has the following important 
fact:

\textit{There is a canonical momentum map one can
  attach to a Poisson action.}

\noindent This is ``the momentum map'' we refer to in the title of this paper.
In order to explain this, we recall \cite{CaFe,CrFe1} that 
\emph{any} Poisson manifold $M$ has an associated symplectic 
groupoid $\Sigma(M)\tto M$, called the \emph{Weinstein groupoid} 
of $M$. In general, the groupoid $\Sigma(M)$ is not smooth but it has the structure of a differentiable symplectic stack (\cite{TsZh1,TsZh2}).
When $\Sigma(M)$ is smooth, $M$ is 
called an \emph{integrable} Poisson manifold and we can think of $\Sigma(M)$ as the
\emph{symplectization} of $M$.
We will see that a Poisson action of $G$ on $M$ always 
lifts to a \emph{Hamiltonian} action of $G$ on $\Sigma(M)$.

\begin{thm}[Symplectization of Poisson actions]
  \label{thm:symp:reduction}
  For a Poisson $G$-space $M$, its symplectic groupoid 
  $\Sigma(M)$ is a Hamiltonian $G$-space that has a natural 
  equivariant momentum map $J:\Sigma(M)\to\gg^*$, which is a groupoid 1-cocycle, that is, 
  \[ J(x\cdot y)=J(x)+ J(y),\quad \text{for any }x,y\in\Sigma(M).\]
\end{thm}

We emphasize that there are no choices involved: $J$ always exists
and is completely natural. Theorem \ref{thm:symp:reduction} is, in essence,
due to Weinstein \emph{et al.} (see \cite{CoDaWe,MiWe,WeXu}). Our point of 
view here is slightly different. We give an explicit simple formula 
for the momentum map $J$ taking advantage of the description of 
$\Sigma(M)$ in terms of cotangent paths which makes sense even in 
the non-smooth case (here one uses the differentiable symplectic stack structure).

Let us assume now that $\Phi:G\times M\to M$ is a \emph{proper and
free} Poisson action, so that the orbit space $M/G$ is also a 
Poisson manifold. If $M$ is an integrable Poisson manifold, then
the quotient $M/G$ is also an integrable Poisson manifold. In this
situation there are two natural groupoids associated with $M/G$:
\begin{enumerate}
\item[(i)] Since $M/G$ is integrable, $\Sigma(M/G)$ is a symplectic
	Lie groupoid integrating the Poisson manifold $M/G$.
\item[(ii)] The lifted action $G\times \Sigma(M)\to\Sigma(M)$ is also
	proper and free. The symplectic quotient:
 	\[ \Sigma(M)/\!/G:=J^{-1}(0)/G.\]
	is a symplectic Lie groupoid that also integrates $M/G$.
\end{enumerate}
It will be shown that these two groupoids have isomorphic Lie algebroids.
A natural question is whether these two groupoids are isomorphic, i.e., 
\begin{equation}
  \label{eq:reduction:commutes} 
  \Sigma(M/G)=\Sigma(M)/\!/G.
\end{equation}
in other words, \emph{does symplectization commute with reduction}?
We will see that, in general, there are topological 
conditions for this equality to hold. In order to describe them,
let us denote by $j:T^*M\to \gg^*$ the momentum map for 
the lifted cotangent action, which is given by
\[ \langle j(\al),\xi\rangle=\al(X_\xi), \quad\al\in T^*M,\; \xi\in\gg,\]
where $X_\xi\in\X(M)$ denotes the infinitesimal generator of the 
action for $\xi\in\gg$.

\begin{thm}[Symplectization commutes with reduction]
  \label{thm:symp:reduct}
  Let $G\times M\to M$ be a proper and free Poisson action. Then 
  symplectization and reduction commute if and only if the
  following groups
  \[ 
  K_{m}:=\frac{\{a:I\to j^{-1}(0) \mid a 
  \text{\rm  ~is a cotangent loop such that }a\sim 0_m\}}
{\{\text{\rm cotangent homotopies with values in } j^{-1}(0)\}}
  \]
  are trivial, for all $m\in M$.
\end{thm}

A simple instance when all $K_m$ vanish is the case of the 
trivial Poisson manifold $M$. In this case we have $\Sigma(M)=T^*M$ and we obtain the 
well-known fact that $T^*(M/G)=T^*M/\!/G$~(\cite{fom}).

In the case of Hamiltonian $G$-spaces the groups $K_m$ have a special simple form, 
since they can be described in terms of the fibers of the momentum map.

\begin{cor}
Let $G\times M\to M$ be a proper and free Hamiltonian action on a symplectic manifold $(M,\omega)$
with momentum map $\mu:M\to\gg^*$. Then symplectization and reduction commute 
if and only if the following groups
\[ K_{m}:=\Ker i_*\subset \pi_1(\mu^{-1}(c),m)\]
are trivial for all $m\in M$, where $c=\mu(m)$ and $i:\mu^{-1}(c)\hookrightarrow M$ is the 
inclusion. 
\end{cor}

For example, these groups vanish if the fibers of the momentum map are simply connected,
or if the second relative homotopy groups of the fibers vanish. The latter occurs when the 
group is compact and the momentum map is proper.

\begin{cor}
Let $G\times M\to M$ be a free Hamiltonian action of a compact Lie group on a 
symplectic manifold $(M,\omega)$ with a proper momentum map $\mu:M\to\gg^*$. Then
symplectization and reduction commute. Moreover, the isotropy groups 
$(\Sigma(M)/\!/G)_{[m]}$
all have the same number of connected components, that is, 
\[ \pi_0\left((\Sigma(M)/\!/G)_{[m]}\right) =\pi_1(M/G)=\pi_1\left( M_{\red},[m]\right), \]
where $M_{\red}=\mu^{-1}(\mathcal{O}_c)/G$ is the symplectic reduced space at value $c=\mu(m)$.
\end{cor}

Note that, in general, we \emph{do not} have $\pi_1(M)=\pi_1(M/G)$, contrary to what happens for 
Hamiltonian actions of compact Lie groups on compact symplectic manifolds (see \cite{Li}).

Let us remove now the assumption that the action is free. For proper actions, 
the quotient $X:=M/G$ is a smooth stratified space
$X=\bigcup_{i\in I} X_i$.
The strata $X_i$ are the connected components of the orbit types $M_{(H)}/G$ 
and the space of smooth functions $C^\infty(M/G)$ can be identified with 
the space $C^\infty(M)^G$ of smooth $G$-invariant functions on $M$. Hence, 
when $G$ acts by Poisson diffeomorphisms, the space of smooth functions 
$C^\infty(M/G)$ becomes a Poisson algebra. The Poisson geometry of $X=M/G$ 
has the following  simple description.,

\begin{thm}[Poisson Stratification Theorem]
  \label{thm:Poisson:reduction}
  Let $G$ act properly by Poisson diffeomorphisms on a Poisson
  manifold $M$. Then the orbit type stratification of $X=M/G$ is a Poisson
  stratification, that is, each orbit type stratum $X_i$ is a smooth Poisson
  manifold and the inclusion $X_i\hookrightarrow X$ is a Poisson map.
\end{thm}

Surprisingly, we could not find in the literature this simple and
clear statement concerning \emph{singular Poisson reduction}, 
which is a generalization of the symplectic stratification
theorem of Sjamaar and Lerman \cite{LeSj} (see, however,
\cite{Sn} where one can find a different approach to singular Poisson reduction). 
In the symplectic case, the key fact in constructing this stratification is a theorem due to 
Guillemin and Sternberg (see \cite{GuSt1}, Theorem 3.5)
stating that the connected components of the \emph{isotropy type
manifolds} $M_H:=\{m \in M\mid G_m=H\}$ are symplectic submanifolds
of $M$. This result has a generalization to Poisson
manifolds where, now, each connected component of $M_H$ is a 
\emph{Lie-Dirac submanifold} (these are the analogs of symplectic 
submanifolds in Poisson geometry; see Appendix \ref{appendix:submanifolds}). 
We emphasize that the inclusion map $M_H\hookrightarrow M$ is \emph{not} a 
Poisson map. From this extension of the Guillemin and Sternberg
result, Theorem \ref{thm:Poisson:reduction} follows in a
straightforward way.

The question of integration (or symplectization) of a Poisson stratified 
space leads naturally to the concepts of \emph{stratified Lie algebroids} 
and \emph{stratified symplectic groupoids}. For example, for a non-free proper action, the Weinstein groupoid $\Sigma(M/G)$ 
of the orbit space $M/G$  is a stratified 
symplectic groupoid. Using the language of stratified Lie theory, we will 
explain how to make sense of the statement ``symplectization commutes with reduction''
in the non-free case.

A natural question that arises in this context is under what conditions two 
different proper Poisson actions lead to the same stratified Poisson spaces.
Loosely speaking, we say that two proper Poisson spaces $(M_1,\{~,~\}_1,G_1)$ and 
$(M_2,\{~,~\}_2,G_2)$ are \emph{Morita equivalent} if the corresponding
action groupoids are Morita equivalent and the equivalence bi-module preserves
Poisson structures. We call a Morita equivalence class a \emph{Poisson orbispace}.
In this context, we have the following result.

\begin{thm}[Poisson orbispaces]
The Poisson stratifications, the stratified Lie algebroids, and 
the stratified symplectic groupoids of Morita equivalent Poisson spaces
are naturally isomorphic.
\end{thm} 

It follows that all these are well defined objects attached to a 
Poisson orbispace. Moreover, we will see that every Poisson orbispace has an underlying Lie 
pseudo-algebra (i.e, an algebraic version of a Lie algebroid), so one can 
even avoid altogether singular stratifications. However, at this point we do 
not know what object integrates this Lie pseudo-algebra.

A very natural issue that will be dealt with in a separate work is the  
convexity property of the momentum map $J:\Sigma(M)\to\gg^*$. For example, 
it is easy to see that the image $J(M)\subset\gg^*$ is a star shaped region which is 
symmetric with respect to the origin. This should lead to nice generalizations of 
the classical Atiyah-Guillemin-Kirwan-Sternberg convexity theorem.
Note also that the results in this paper can be extended in various 
directions by replacing Poisson structures by Dirac structures, Jacobi structures,
and other higher geometric structures.

One can also place our work in a wider context by allowing general Poisson actions 
by Poisson-Lie groups in the sense of Jiang-Hua Lu \cite{Lu}. Our results do extend to 
this more general setting. Indeed, in this case, the map $j:T^*M \to\gg^*$ above is a 
morphism from the Lie bialgebroid $(T^*M,TM)$ to the Lie bialgebra $(\gg^*,\gg)$. Assuming
that $(M,\Pi)$ is an integrable Poisson manifold, $j$ integrates to a morphism $J:\Sigma(M)\to G^*$
from the symplectic groupoid of $M$ to the dual Poisson-Lie group $G^*$, and the lifted action
is a Poisson-Lie group action of $G$ on $\Sigma(M)$ with momentum map $J$. The reduced 
symplectic groupoid integrating $M/G$ is now $\Sigma(M)/\!/G:=J^{-1}(e)/G$. This generalization
requires techniques related to double structures which are quite different from those used in this paper.  
We refer to the upcoming works \cite{FeIg,Stef} for details.

%\vskip 5 pt

The paper is organized as follows. Singular Poisson reduction and the proof of
the Poisson Stratification Theorem are presented in Section \ref{sec:singular:reduction}. The symplectization of Poisson actions
and of the associated momentum maps is studied in Section
\ref{sec:momentum:maps}. The problem of ``symplectization commutes with reduction'' is discussed in detail in Section \ref{sec:symp:reduction}. The last section is  dedicated to the study of Poisson orbispaces.

\vskip 5 pt

\noindent
\textbf{Acknowledgments.} We thank the referee for useful remarks and suggestions and Marius Crainic, Markus Pflaum, Miguel Rodr\'{\i}guez-Olmos, Alan Weinstein, and Nguyen Tien Zung for enlightening discussions. The authors thank the Bernoulli Center and the Mathematics Department of IST for its hospitality, where parts of this work were done. 

\tableofcontents

%%%%%%%%%%%%%%%%%%%%%%%%%%%%%%%%%%%%
%%%%%%%%%%%%%%%%%%%%%%%%%%%%%%%%%%%%
\section{Singular Poisson reduction}%
\label{sec:singular:reduction}     %
%%%%%%%%%%%%%%%%%%%%%%%%%%%%%%%%%%%%
%%%%%%%%%%%%%%%%%%%%%%%%%%%%%%%%%%%%
%%%%%%%%%%%%%%%%%%%%%%%%%%%%%%%%%%%%

The purpose of this section is to study the structure of the orbit
space $M/G $ when $(M,\{\cdot,\cdot\})$ is a Poisson manifold and
$G$ is a Lie group that acts properly by Poisson diffeomorphisms 
on $(M,\{\cdot,\cdot\})$. 

Recall that for proper actions we have a classical result 
(see, e.g., Theorem 2.7.4 in  \cite{DuKo} as well as the discussion following that result) stating that the connected 
components of the \emph{orbit type} manifolds $M_{(H)}$ and their 
projections $M_{(H)}/G$, constitute Whitney stratifications of $M$ and 
$M/G$, respectively, as $H$ varies in the set of all 
isotropy subgroups of $G$. Our aim is to show that the Poisson 
geometry of the smooth stratified space $M/G$ has a simple description.

We start by proving a general result about the Poisson
nature of the fixed point set of a Poisson action. This result leads
to the Poisson version of a theorem, due to Guillemin and Sternberg
(see \cite[Theorem 3.5]{GuSt1}), that shows that the connected
components of the \emph{isotropy type manifolds} 
$M_H:=\{m \in M\mid G_m=H\}$ are symplectic submanifolds of $M$. 
Using the result on the Poisson nature of the fixed point sets 
one constructs Poisson structures for the orbit type strata of $M/G$. 
In the last paragraph of the section, we will describe the 
symplectic leaves of these strata using the 
\emph{optimal momentum map} introduced in \cite{OrRa1}.

%%%%%%%%%%%%%%%%%%%%%%%%%%%%%%%%%%%%%%%%%%%%%%%%%
\subsection{The fixed point set of a Poisson action}%
\label{subsec:fixed:points}                     %
%%%%%%%%%%%%%%%%%%%%%%%%%%%%%%%%%%%%%%%%%%%%%%%%%

Let $\Phi:G\times M\to M$ be a smooth, proper, Poisson action on the
Poisson manifold $(M, \{\cdot,\cdot\})$. Let $\Pi\in\X^2(M)$ be the
associated Poisson tensor given by $\Pi(\d f ,\d h):=\{f,h\}$, 
for any $f,h\in C^\infty(M)$, and $\sharp:T^\ast M\rightarrow TM$ 
the vector bundle map defined by  
\[
\sharp (\d  f )=X_f:=\{ \cdot , f \}.
\]
We denote by $M^G:=\{m\in M\mid g\cdot m=m,\text{ for all } g\in G\}$
the \emph{fixed point set} of the action. Since the action is proper,
each connected component of $M^G$ is a submanifold of $M$. Actually,
if $M^G$ is non-empty, then $G$ must be compact, so we will assume
compactness in the following result (for the terminology used in the
statement we refer to Appendix~\ref{appendix:submanifolds}):

\begin{thm}
  \label{Lie-Dirac fixed point spaces}
  Let $G$ be a compact Lie group and $M$ a Poisson $G$-space. Then $M^G$ is a
  Lie-Dirac submanifold of $M$ with Poisson bracket $\{\cdot,\cdot\}_{M^G}$ 
  given by
  \begin{equation}
    \label{Poisson bracket of fixed point manifold}
    \{f,h\}_{M^G}:=\left.\{\widetilde{f},\widetilde{h}\}\right|_{M^G},\quad
    f,h\in C^\infty(M^{G}),
  \end{equation}
  where $\widetilde{f},\widetilde{h}\in C^\infty(M)^{G}$ denote arbitrary
  $G$-invariant extensions of $f,h\in C^\infty(M^{G})$.
\end{thm}

For the proof, we fix, once and for all, a $G$-invariant Riemannian metric
$(~,~)$ on $M$. Let  
\[
E=\{v\in T_{M^G}M \mid (v,w)=0,\forall w\in TM^G\} \subset T_{M^G}M
\]
be the orthogonal subbundle to $TM^G$.

\begin{lem}
$T_{M^G} M=TM^G\oplus E \text{ and } \;\sharp(E^0)\subset TM^G.$
\end{lem}

\begin{proof}
Since $E=(TM^G)^\perp$, the decomposition 
$T_{M^G}M=TM^G\oplus E$
is obvious. Moreover, since $G$ is compact, the action
linearizes around fixed points and we see that, for the lifted tangent
action, $(T_{M^G}M)^G=TM^G$. It follows that this decomposition can also be written as
\begin{equation}
\label{eq:Decomp1}
T_{M^G} M=(T_{M^G}M)^G\oplus E.
\end{equation}
On the other hand, the lifted cotangent action $G\times
T^*M\to T^*M$ is related to the lifted tangent action by
$g\cdot\xi(v)=\xi(g^{-1}\cdot v)$, $\xi\in T^*M, v\in TM$.
We claim that
\begin{equation}
\label{eq:Decomp2} 
E^0=(T^*_{M^G}M)^G,
\end{equation}
where $E^0 $ is the annihilator of $E $ in $T^*_{M^G}M $.
Indeed, if $v\in T_{M^G}M$ we can use (\ref{eq:Decomp1}) to decompose
it as $v=v_G+v_E$, where $v_G\in(T_{M^G}M)^G$ and $v_E\in E$. Hence, for $\xi\in E^0$ we find 
\begin{align*}
g\cdot\xi(v_G+v_E)
&=\xi(g^{-1}\cdot v_G+g^{-1}\cdot v_E)
=\xi(v_G)+\xi(g^{-1}\cdot v_E) 
=\xi(v_G)\\
&=\xi(v_G)+\xi(v_E)=\xi(v_G+v_E).
\end{align*}
We conclude that $g\cdot\xi=\xi$ for all $\xi\in E^0$ and hence $E^0\subset (T^*_{M^G}M)^G$. Now
(\ref{eq:Decomp2}) follows by counting fiber dimensions.

Since $G\times M\to M$ is a Poisson action, $\sharp:T^*M\to TM$ is a
$G$-equivariant bundle map. Hence, if $\xi\in E^0$, we obtain 
from (\ref{eq:Decomp2}) that
$g\cdot\sharp\xi=\sharp(g\cdot\xi)=\sharp\xi$.
This means that $\sharp\xi\in (T_{M^G}M)^G=TM^G$, so the lemma holds.
\end{proof}

\begin{rem}
\label{rem:metric}
Equations (\ref{eq:Decomp1}) and (\ref{eq:Decomp2}) show that $E=\left[(T^*_{M^G}M)^G\right] ^0$. Hence, even though we used a metric to introduce $E$, this bundle does not depend on the particular choice of metric. This also shows that
the way in which $M^G$ embeds in $M$ as a Poisson-Dirac submanifold is also independent of the choice of metric.
\end{rem}

The previous lemma shows that the conditions of Corollary \ref{cor:LieDirac}
in the appendix are satisfied, so $M^G$ is a Poisson-Dirac submanifold. 
Formula (\ref{Poisson bracket of fixed point manifold}) for
the bracket follows from equation (\ref{eq:PoissonBracket}) in the
appendix, together with the fact that for any $G$-invariant extension
$\widetilde{f}$ of $f\in C^\infty(M^{G})$ we have 
$\d_x\widetilde{f}\in E^0$, for $x\in M^G$. Actually, the sections of 
$E^0$ are generated by the differentials $\d_{M^G}\widetilde{f}$ of 
$G$-invariant extensions of functions $f\in C^\infty(M^{G})$. We have the following result. 

\begin{lem}
\label{E_zero_algebroid}
  $E^0$ is a Lie subalgebroid of $T^*M$.
\end{lem}

\begin{proof}
The canonical identification of $E^0$ with $T^*M^G$
defines the Lie algebroid structure on $M^G$, namely, the cotangent Lie
algebroid structure of the Poisson structure on $M^G$ (see, also, Remark
\ref{rem:metric} above). Therefore, we have to show that the inclusion 
$T^*M^G\simeq E^0\hookrightarrow T^*M$ is a Lie algebroid morphism.

We already know that the anchor $\sharp_M$ maps $E^0$ into $TM^G$. Since,
under the identification $T^*M^G\simeq E^0$, we have $\sharp_{M^G}=\sharp_M|_{E^0}$,
we conclude that inclusion preserves anchors. It remains to show that
the inclusion preserves brackets. To explain what this means, given any
two sections $\al,\be\in\Gamma(T^*M^G)\simeq\Gamma(E^0)$, we can
always write them as 
\[ \al=\sum_i a_i\, i^*\al_i,\quad \al=\sum_j b_j\, i^*\be_j,\]
where $a_i,b_j\in C^\infty(M^G)$, $\al_i,\be_j\in\Omega^1(M)$, and
$i:M^G\to M$ is the inclusion. The requirement for the inclusion
to preserve brackets is the following identity:
\[ [\al,\be]_{T^*M^G}=\sum_{ij}a_i b_j i^*[\al_i,\be_j]_{T^*M}
+\sum_j\sharp\al(b_j) \be_j-\sum_i\sharp\be(a_i) \al_i\]
(this is independent of the choices of $a_i,\al_i,b_j,\be_j$).
Let us show that this equality holds. Note that, for any section
$\al\in\Gamma(T^*M^G)$, we can write
\[ \al=\sum_i a_i \d f_i=\sum_i a_i\, \d_{M^G}\widetilde{f_i}
=\sum_i a_i\, i^*\d \widetilde{f_i},\]
where $a_i,f_i\in C^\infty(M^{G})$ and $\widetilde{f_i}\in
C^\infty(M)$ are $G$-invariant extensions. Similarly, if
$\be\in\Gamma(T^*M^G)$ is another section, we have 
$\be=\sum_j b_j\, i^*\d \widetilde{g_j}$,
for $b_j,g_j\in C^\infty(M^{G})$, so we get
\begin{align*}
[\al,\be]_{T^*M^G}
 &=\sum_{ij}{\left[a_i\d f_i,b_j\d g_j\right]}_{T^*M^G}\\
 &=\sum_{ij}a_i b_j[\d f_i,\d g_j]_{T^*M^G}+
    \sum_j\sharp\al(b_j)\d g_j-\sum_i\sharp\be(a_i)\d f_i\\
 &=\sum_{ij}a_i b_j\d\{f_i, g_j\}_{M^G}+
    \sum_j\sharp\al(b_j)\d_{M^G}\widetilde{g_j}-
    \sum_i\sharp\be(a_i)\d_{M^G}\widetilde{f_i}\\
 &=\sum_{ij}a_i b_j\d_{M^G}\{\widetilde{f_i},\widetilde{g_j}\}_M+
    \sum_j\sharp\al(b_j)\d_{M^G}\widetilde{g_j}-
    \sum_i\sharp\be(a_i)\d_{M^G}\widetilde{f_i}\\
 &=\sum_{ij}a_i b_ji^*[\d \widetilde{f_i},\d \widetilde{g_j}]_{T^*M}+
    \sum_j\sharp\al(b_j)i^*\d\widetilde{g_j}-
    \sum_i\sharp\be(a_i)i^*\d\widetilde{f_i}\\
 &=\sum_{ij}a_i b_j i^*[\al_i,\be_j]_{T^*M}+
    \sum_j\sharp\al(b_j) \be_j-\sum_i\sharp\be(a_i) \al_i,
\end{align*}
as required.
\end{proof}

By Lemma \ref{E_zero_algebroid}, $M^G$ is a Lie-Dirac submanifold and this completes 
the proof of Theorem \ref{Lie-Dirac fixed point spaces}.

\begin{rem}
  Special cases of Theorem \ref{Lie-Dirac fixed point spaces} where
  known before. Damianou and Fernandes in \cite{DaFe} show that the
  fixed point set is a Poisson-Dirac submanifold, but we will need the
  stronger statement that $M^G$ is a Lie-Dirac submanifold. Fernandes and
  Vanhaecke  consider in \cite{FeVan} the case where $G$ is a reductive
  algebraic group. Xu in \cite{Xu} proves the Poisson involution
  theorem, which amounts to the case $G=\Zz_2$. Xiang Tang's PhD thesis \cite{Tang} also contains
  a version of this theorem.
\end{rem}

It should be noted that the inclusion $M^G\hookrightarrow M$ is
\emph{not} a Poisson map. It is, in fact, a (backwards) Dirac
map. This means that the graph of the Poisson tensor $\Pi_{M^G}$ is
the pull back of the graph of the Poisson tensor $\Pi_M$:
\[ \text{Graph}(\Pi_{M^G})=
\set{(v,\xi|_{TM^G})\in TM^G\oplus T^*M^G~|~\xi\in T^*M\text{ and }
  v=\Pi_M(\xi,\cdot)}.
\]

We end this paragraph with a simple example.

\begin{example}
Let $\tau:\gg\to\gg$ be an involutive Lie algebra automorphism. Then
$\gg=\hh\oplus\pp$, 
where $\hh$ and $\pp$ are the $\pm 1$-eigenspaces of $\tau$. On $M=\gg^*$ we 
consider the Lie-Poisson structure and we let $\Zz_2=\{1,-1\}$ act on $M$
by $(-1)\cdot\xi\equiv\tau^*(\xi)$. Since $\tau:\gg\to\gg$ is a Lie algebra
automorphism, its transpose $\tau^*:\gg^*\to\gg^*$ is a Poisson
diffeomorphism, so this defines indeed a Poisson action. We conclude that
the fixed point set
\[ M^{\Zz_2}=\pp^0\simeq \hh^*\subset\gg^*,\]
is a Lie-Dirac submanifold. This fact is well known to
people working in integrable systems and is part of the so-called
Adler-Kostant-Symes scheme. 

More generally, we can consider a Lie algebra automorphism
$\tau:\gg\to\gg$ of order $q$. Its fixed point set is a Lie subalgebra
$\hh\subset \gg$. On the other hand, by transposition we obtain a
Poisson action of $\Zz_q$ on $\gg^*$. 
By Theorem \ref{Lie-Dirac fixed point spaces}, its fixed point set 
\[ M^{\Zz_q}=\{\al\in\gg^*:\tau^*(\al)=\al\}\simeq\hh^*\subset\gg^*\]
is a Lie-Dirac submanifold. For a specific example, we can take the
rank 4 orthogonal Lie algebra $\mathfrak{d}_4$ and let
$\tau:\mathfrak{d}_4\to \mathfrak{d}_4$ be the automorphism of order 3 
induced from the obvious $\Zz_3$-symmetry of its Dynkin diagram
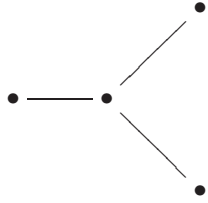
\[
\xymatrix{
&&\bullet\\
\bullet\ar@{-}[r]&\bullet\ar@{-}[ur]\ar@{-}[dr]\\
&&\bullet
}\]
The fixed points of the corresponding $\Zz_3$-action is the
exceptional Lie algebra $\mathfrak{g}_2$. It follows that
$\mathfrak{g}_2^*$ is a Lie-Poisson subspace of $\mathfrak{d}_4^*$.
\end{example}

%%%%%%%%%%%%%%%%%%%%%%%%%%%%%%%%%%%%%%%%%%%%%%%%%%%%%%%%%%%%%
\subsection{Poisson geometry of the isotropy type manifolds}%
\label{subsec:isotropy}                                     %
%%%%%%%%%%%%%%%%%%%%%%%%%%%%%%%%%%%%%%%%%%%%%%%%%%%%%%%%%%%%%

As before, we let $\Phi:G\times M\to M$ be a smooth, proper, Poisson
action. We denote by $G_m$ the \emph{isotropy group} of a point $m\in
M$, by $M_H:=\set{m\in M\mid G_m=H}$ the $H $-\emph{isotropy type submanifold}, 
and by $M^H:=\{m\in M\mid g\cdot m=m,\text{ for all }g\in H\}$ the
$H$-\emph{fixed point manifold}. For a subgroup $H\subset G$ we will
denote by $(H)$ its conjugacy class, and we let 
$M_{(H)}=\{m\in M\mid G_m\in (H)\}$ denote the $(H)$-\emph{orbit type
  submanifold}. The properness of the action
guarantees that each $G_m$ is a compact Lie group and that the
connected components of $M_H$, $M^H$, and $M_{(H)}$ are embedded submanifolds of
$M$. We recall that $M_H$ is an open subset of $M^H$ and that
\[ M_H=M_{(H)}\cap M^H.\]

\begin{rem}
The connected components of $M_H$, $M^H$, and $M_{(H)}$ may be submanifolds
of different dimensions. Henceforth, we will allow our submanifolds
to have components of different dimension (these are sometimes called $\Sigma$-manifolds).
\end{rem}

The following result is the analogue in Poisson geometry of a
well-known theorem in symplectic geometry due to Guillemin and Sternberg
(\cite{GuSt1}, Theorem 3.5).

\begin{prop}
  \label{prop:isotropy:types}
  Let $\Phi:G\times M\to M$ be a proper Poisson action, let
  $H\subset G$ be an isotropy group, and denote by $N(H)$ the normalizer
  of $H$ in $G$. Then:
  \begin{enumerate}[(i)]
  \item $M_H$ is a Lie-Dirac submanifold of $M$ with Poisson
    bracket given by 
    \[
    \{f,h\}_{M_H}=\{\widetilde{f},\widetilde{h}\}|_{M_H},\qquad 
    f,g\in C^{\infty}(M_H),
    \]
    where $\widetilde{f},\widetilde{h}\in C^\infty(M)^H $ denote any
    $H$-invariant extensions of $f,h\in C^{\infty}(M_H)$.
  \item The natural action of $L(H):=N(H)/H$ on $M_H$ is a proper and
    free Poisson action.
  \end{enumerate}
\end{prop}

\begin{proof}
Part (i) follows from Theorem~\ref{Lie-Dirac fixed point spaces} by
replacing the group $G$ by the subgroup $H$. For part (ii) we
recall (see Proposition \ref{prop:Poisson:Dirac:maps} in
Appendix~\ref{appendix:submanifolds}) that a Poisson diffeomorphism
$\phi:M\to M$ leaving a Poisson-Dirac submanifold $N\subset M$ invariant
determines, by restriction, a Poisson diffeomorphism $\phi:N\to N$. Hence, the
action of $N(H)$ on $M_H$ is Poisson, and (ii) follows.
\end{proof}

By a standard result on proper and free Poisson actions, the Poisson
bracket $\{\cdot,\cdot\}_{M_H}$ induces a unique Poisson bracket
$\{\cdot,\cdot \}_{M_H/L(H)}$ on the orbit space $M_H/L(H)$ that
makes the projection $M_H\rightarrow M_H/L(H)$ into a Poisson
map. This will be used in the next paragraph to prove the Poisson
stratification theorem.

%%%%%%%%%%%%%%%%%%%%%%%%%%%%%%%%%%%%%%%%%%%%%%%%%
\subsection{Poisson stratifications}            %
\label{subsec:Poisson:stratifications}          %
%%%%%%%%%%%%%%%%%%%%%%%%%%%%%%%%%%%%%%%%%%%%%%%%%

Recall that if a Lie group $G$ acts properly on a manifold $M$, the 
orbit space $X:=M/G$ is a smooth stratified space (see, e.g, \cite{DuKo,Pflaum}). The decomposition $M=\bigcup_{(H)}M_{(H)}$ into orbit types induces the decomposition 
\[ X=\bigcup_{(H)}M_{(H)}/G\]
of the orbit space $X=M/G$.
The smooth stratification of $X$ is then
\[ X=\bigcup_{i\in I} X_i,\]
where each $X_i$ is a connected component of some $M_{(H)}/G$. The algebra of smooth 
functions on the orbit space $X$ is
\[
C^\infty(X)=\{f\in C^0(M/G)\mid f\circ\pi_G\in C^\infty(M)^{G}\}.
\]
We will show that this is a \emph{Poisson stratification} in the sense
of the following definition.

\begin{definition}
  Let $X$ be a topological space. A \textbf{Poisson stratification} of
  $X$ is a smooth stratification $\mathcal{S}=\{S_i\}_{i\in  I}$ of $X$ together with a
  Poisson algebra $(C^{\infty}(X),\{\cdot,\cdot\}_X)$, where 
  $C^{\infty}(X)\subset C^0(X)$ is the space of smooth functions associated with $\mathcal{S}$, 
  such that:
  \begin{enumerate}[(i)]
  \item Each stratum $S_i$, $i\in I$, is a Poisson manifold. 
  \item The inclusions $i:S_i \hookrightarrow X$ are
    Poisson maps, i.e., $\{f,h\}_X \circ i=\{f \circ i,h\circ i\}_{S_i}$, 
    for any  $f,h\in C^{\infty}(X)$ and $i\in I$. 
  \end{enumerate}
\end{definition}

\begin{rem}
  \label{rem:uniqueness}
  Note that, given a smooth stratification $\{S_i\}_{i\in I}$ of a
  topological space $X$ together with a Poisson bracket on its algebra
  of smooth functions $C^{\infty}(X)$, there is at most one structure of
  a Poisson stratification on $\{S_i\}_{i\in I}$. In other words, the
  Poisson structures on the strata $S_i$ are uniquely determined by the
  Poisson algebra $(C^{\infty}(X),\{\cdot,\cdot\}_X)$.
\end{rem}

For a proper Poisson action $\Phi:G\times M\to M$ the smooth functions
on $M/G$ have a natural Poisson algebra structure, namely the one
obtained by identifying $C^\infty(M/G)\simeq C^\infty(M)^{G}$ (it is
the unique one for which the natural projection $\pi_G:M\rightarrow
M/G$ is a Poisson map). Let us show now that the strata of $M/G$ have
canonical Poisson structures.

\begin{prop}
  \label{strata Poisson}
  Let $\Phi:G\times M\to M$ be a proper Poisson action and $H\subset
  G$ an isotropy group.
  \begin{enumerate}[(i)]
  \item The natural map $F_H:M_H/L(H)\rightarrow M_{(H)}/G$ is a
    diffeomorphism, so $M_{(H)}/G$ inherits a Poisson structure from
    $M_H/L(H)$. 
  \item If $H_1,H_2\in(H)$ are conjugate isotropy groups, the Poisson
     structures on $M_{(H)}/G$ induced from $M_{H_1}/L(H_1)$ and
     $M_{H_2}/L(H_2)$ coincide.
  \end{enumerate}
\end{prop}

\begin{proof}
The fact that the map $F_H:M_H/L(H)\rightarrow M_{(H)}/G$ is a
diffeomorphism is standard. To prove (ii), we show that the 
isomorphism $\phi:M_{H_1}/L(H_1)\to M_{H_2}/L(H_2)$ obtained by 
composition
\[ 
\xymatrix{
M_{H_1}/L(H_1)\ar[rr]^{F_{H_1}}& &M_{(H)}/G \ar[rr]^{F_{H_2}^{-1}}&& M_{H_2}/L(H_2)
}
\]
is Poisson. To see this, let $g\in G$ be such that $H_2=gH_1 g^{-1}$
and let $\Phi_g:M\to M$ be defined by $\Phi_g(x)=g\cdot x$, for any $x \in M $.
Then we have a commutative diagram
\[
\xymatrix{ 
M_{H_1}\ar[d]\ar[r]^{\Phi_g}&M_{H_2}\ar[d]\\
M_{H_1}/L(H_1)\ar[r]_{\phi}& M_{H_2}/L(H_2)}
\]
where the top row is a Poisson diffeomorphism and the vertical
projections are Poisson. Hence, $\phi$ must be a Poisson
diffeomorphism as well.
\end{proof}

The Poisson structure on each stratum can also be obtained by Dirac
reduction. Namely, each orbit type submanifold $M_{(H)}$ gets an induced
Dirac structure 
\[ L_{M_{(H)}}:=\{(v,\eta)\in TM_{(H)}\oplus T^*M_{(H)}\mid \exists \xi\in T^*M
\text{ such that }\xi|_{TM_{(H)}}=\eta,~\sharp\xi=v\}\]
from the Poisson structure on $M$.
Note that this is just the pull back Dirac structure $i^*L_\pi$, where
$L_\pi=\text{graph}(\pi)$ and $i:M_{(H)}\hookrightarrow M$. Now, this
pull back Dirac structure is clearly $G$-invariant and descends 
to the quotient $M_{(H)}/G$. It is easy to check that the reduced 
Dirac structure is just the graph of the reduced Poisson structure on $M_{(H)}/G$
we have constructed above. 

All this can be expressed by the commutative diagram of Dirac manifolds
\[
\xymatrix{
     &M\\
M_H\ar[ur]\ar[dr]& &M_{(H)}\ar[ul]\ar[dl]\\
 &M_H/L(H)\simeq M_{(H)}/G}
\]
where the inclusions are backward Dirac maps and the projections are
forward Dirac maps.

Now we can prove the Poisson Stratification Theorem.

\begin{thm}
  \label{Poisson strata we go}
  Let $\Phi:G\times M\to M$ be a proper Poisson action. The connected
  components of the orbit type reduced spaces $M_{(H)}/G$ form a 
  Poisson stratification of $(M/G,\{\cdot,\cdot\}_{M/G})$.
\end{thm}

\begin{proof}
The classical Orbit Type Stratification Theorem proves the 
stratification part of the statement. It remains to be shown that 
the inclusions $i:M_{(H)}/G\hookrightarrow M/G $ are Poisson 
maps, i.e., that 
\[ \{f,h\}_{M/G}\circ i=\{f\circ i, h\circ i\}_{M_{(H)}/G},\]
for any $f,h\in C^{\infty}(M/G)$.

To see this, we consider, as before, the isomorphism 
$F_H:M_H/L(H)\to M_{(H)}/G$ and let $\pi_{L}:M_H\rightarrow M_H/L(H)$
be the projection. Then, for any $m\in M_H$, we have
\begin{align*}
\{f\circ i, h\circ i\}_{M_{(H)}/G}(F_H([m]))
&=\{f\circ i\circ F_H, h\circ i\circ F_H\}_{M_H/L(H)}(\pi_L(m))\\
&=\{f\circ i\circ F_H\circ\pi_L, h\circ i\circ F_H\circ\pi_L\}_{M_H}(m)\\
&=\{f\circ\pi_G,h \circ\pi_G\}_M(m)\\
&=\{f,h\}_{M/G} \left(\pi_G(m) \right) = \left( \{f,h\}_{M/G} \circ i \right) (F_H([m]),
\end{align*}
where we have used the fact that 
$f\circ\pi_G,h\circ\pi_G\in C^\infty(M)$ are $G$-invariant (and hence
$H$-invariant) extensions of $f\circ i\circ F_H\circ \pi _L,
h \circ i\circ F_H\circ\pi_L\in C^{\infty}(M_H)$, respectively.
\end{proof}

\begin{remark}
The Poisson structure of the strata in the previous theorem can also
be obtained by using a combination of results that one can find in the
literature. First, \cite[Proposition 5]{Sn} proves in the context of
subcartesian Poisson spaces that the accessible sets by the Hamiltonian flows
in $(M/G, C^\infty(M/G))$ carry a natural Poisson structure. Second, the main theorem 
in \cite{Or} (see also \cite[Theorem 3.5.1]{OrRa}) proves that those
accessible sets are precisely the connected components of the orbit
spaces $M_H/L(H)$.
\end{remark}

%%%%%%%%%%%%%%%%%%%%%%%%%%%%%%%%%%%%%%%%%%%%%%%%%
\subsection{Symplectic leaves of the strata}    %
\label{subsec:optimal:momentum}                 %
%%%%%%%%%%%%%%%%%%%%%%%%%%%%%%%%%%%%%%%%%%%%%%%%%

We describe now the symplectic leaves
of the Poisson strata $M_H/G$ of $(M/G,\{\cdot,\cdot\}_{M/G})$
introduced in the Poisson Stratification Theorem. In order to achieve
this we will make use of the \emph{optimal momentum map} and the
\emph{optimal reduction}, introduced in \cite{OrRa1} and~\cite{Or1}, respectively, which we briefly
review.

As before, we let $\Phi: G\times M \rightarrow M$ be a proper Poisson
action. Let  $A_{G}:=\{\Phi_g\mid  g \in G \}$ be the
associated group of Poisson diffeomorphisms and $A'_G $ the integrable
generalized distribution defined by
\[ A_{G}':=\left\{X _f\mid f \in C^\infty(M)^G\right\}.\] 
The \textbf{optimal momentum map} $\mathcal{J}:M\rightarrow  M/A_{G}'$
of this Poisson action is defined as the projection of $M$ onto the
leaf space $M/A_{G}'$, endowed with the quotient topology (note that
this space can be quite singular!). The main facts concerning 
\textbf{optimal reduction} are the following (see \cite{OrRa1, Or1}):
\begin{enumerate}[(i)]
\item There exists a unique continuous $G$-action on $M/A'_G$ that
  makes the optimal momentum map $G$-equivariant. 
\item For any $\rho \in M / A_{G}'$, the isotropy subgroup $G _\rho$
  has a unique smooth structure that makes it into an initial
  submanifold of $G$ (recall that a submanifold $i:N\hookrightarrow M $ is \emph{initial} if
  the inclusion is a regular immersion, i.e., for every smooth manifold $P $, a map
  $f:P\rightarrow N $ is smooth if and only if $i \circ  f :P \rightarrow M $ is smooth). 
\item For any $\rho \in M / A_{G}'$, $\mathcal{J}^{-1}(\rho)$ is an
  initial submanifold of $M$.
\item If $G_{\rho}$ acts properly on $\mathcal{J}^{-1}(\rho)$ the
  orbit space $M_\rho:=\mathcal{J}^{-1}(\rho)/G_\rho$ is a smooth
  symplectic manifold with symplectic form $\omega_\rho$ defined by:
  \begin{equation*}
    \label{symplectic 1}
    \qquad(\pi_\rho^\ast\omega_\rho)(m)(X_f(m),X_h(m))=\{f,h\}_M(m), \quad
    (m\in \mathcal{J}^{-1}(\rho),\, f,h\in C^\infty(M)^G)
  \end{equation*} 
  where $\pi_\rho:\mathcal{J}^{-1}(\rho)\to M_\rho $ is the canonical
  projection. The pair $(M_\rho,\omega_\rho)$ is called the
  \textbf{optimal point reduced space} at $\rho$.
\item If  $\mathcal{O}_{\rho}=G\cdot\rho\subset M/A_{G}'$ is the
  $G$-orbit of $\rho\in M/A_{G}'$, the map 
  \[ 
  \mathcal{J}^{-1}(\rho)/G_{\rho}\to\mathcal{J}^{-1}(\mathcal{O}_{\rho})/G, 
  \quad [m]_\rho\longmapsto [m]_{\mathcal{O}_\rho},
  \] 
  is a bijection, so the quotient
  $M_{\mathcal{O}_\rho}:=\mathcal{J}^{-1}(\mathcal{O}_{\rho})/G$ has a
  smooth symplectic structure $\omega_{\mathcal{O}_\rho}$ induced from
  $(M_\rho,\omega_\rho)$. The pair
  $(M_{\mathcal{O}_\rho},\omega_{\mathcal{O}_\rho})$ is called the 
  \textbf{optimal orbit reduced space} at $\rho$.
\end{enumerate}

The symplectic foliation of the Poisson stratified space $M/G$ can
now be described as follows.

\begin{thm}
\label{thm:symplectic:leaves}
Let $\Phi: G\times M \rightarrow M$ be a proper Poisson action with
optimal momentum map $\mathcal{J}:M\rightarrow M/A_{G}'$. The
symplectic leaf of the stratum $M_{(H)}/G$ through $[m]$ is the optimal
orbit reduced space 
$\left(\mathcal{J}^{-1}(\mathcal{O}_{\rho})/G,\omega_{\mathcal{O}_\rho}\right)$
at $\rho=\mathcal{J}(m)$.
\end{thm}

For the proof we need the following two lemmas.

\begin{lem}
The optimal orbit reduced spaces are contained in the strata:
\[
M_{\mathcal{O}_\rho}=\mathcal{J}^{-1}(\mathcal{O}_{\rho})/G\subset M_{(H)}/G. 
\] 
\end{lem}

\begin{proof}
The equivariance of $\mathcal{J}$ implies that 
$\mathcal{J}^{-1}(\rho)\subset M_H$, so we conclude that
\[
\mathcal{J}^{-1}(\mathcal{O}_{\rho})/G=G\cdot\mathcal{J}^{-1}(\rho)/G 
\subset G \cdot M_H/G=M_{(H)}/G
\] 
and the lemma follows.
\end{proof}

\begin{lem}
Let $\leaf_{[m]}$ be the symplectic leaf of $M_{(H)}/G$ through $[m]$. Then
\[ T_{[m]}\leaf_{[m]}=T_{[m]}M_{\mathcal{O}_\rho}.\]
\end{lem}

\begin{proof}
Let  $m\in M_H$ and $f\in C^{\infty}(M_H)^{N(H)}$.
Proposition 2.5.6 in \cite{OrRa} guarantees the existence of an open
$G$-invariant neighborhood $U$ of $m$ in $M$ and of a $G$-invariant
function $\widetilde{f}\in C^{\infty}(U)^G $  such that
$\widetilde{f}|_{U\cap M_H}=f|_{U\cap M_H}$. Moreover, if $X_f$
and $X_{\widetilde{f}}$ denote the Hamiltonian vector fields
associated to $f$ and $\widetilde{f}$ with respect to 
$\{\cdot,\cdot\}_{M_H}$ and $\{\cdot,\cdot\}_M$, they necessarily
coincide on $U\cap M_H$, that is, 
\begin{equation}
\label{equality of vector fields}
X_f=X_{\widetilde{f}}.
\end{equation}
Indeed, due to the $G$-invariance of $\widetilde{f}$, the vector field 
$X_{\widetilde{f}}$ is tangent to $M_H$ when evaluated at points in
$U\cap M_H$. Since $M_H$ is a Poisson-Dirac submanifold of $M$, its symplectic
leaves are the (connected components of the) intersection of $M_H$
with the symplectic leaves of $M$, so 
(\ref{equality of vector fields}) follows.

Let $\pi_G:\mathcal{J}^{-1}(\mathcal{O}_{\rho}) \rightarrow M
_{\mathcal{O}_\rho}:=\mathcal{J}^{-1}(\mathcal{O}_{\rho})/G$ be the
projection. The tangent space of $M_{\mathcal{O}_\rho} $ at any point
$\pi_G (m)\in M_{\mathcal{O}_\rho}$ is the set of vectors of the form
$T_m\pi_G \cdot X_l(m)$, where $l\in C^\infty(M)^{G}$ is
arbitrary. If $m\in M_H$, then we can further assume that
$l=\widetilde{f}$, a $G$-invariant extension of some smooth function
$f\in C^{\infty}(M_H)^{N(H)}$. Hence,
\begin{equation}
\label{tangent:1}
T_{[m]} M_{\mathcal{O}_\rho}=\set{ T_m\pi_G \cdot  X_{\widetilde{f}}(m) \mid  f\in C^{\infty}(M_H)^{N(H)} }.
\end{equation}

On the other hand, the characteristic distribution of the quotient
Poisson manifold $M_H/L(H)$ at $\pi_{L}(m)$ consists of the vectors of
the form $T_m\pi_{L ^m}\cdot X_f(m)$, with $f\in
C^{\infty}(M_H)^{N(H)}$, which by (\ref{equality of vector fields})
equals $T_m\pi_{L} \cdot X_{\widetilde{f}}(m)$, $ \widetilde{f}\in
C^\infty(M)^{G}$.  Consequently, if $F_H$ is the diffeomorphism in
Proposition~\ref{strata Poisson}, the tangent space to the symplectic
leaf of $M_{(H)}/G$ at $F_H([m])$ is
\begin{equation}
\label{tangent:2}
T_{F_H([m])}\leaf=
\set{T_{\pi_{L}(m)}F_H \cdot T_m \pi_{L} \cdot X_{\widetilde{f}}(m) \mid
 \widetilde{f} \in C^\infty(M)^{G}}.
\end{equation} 
Since we have $T_{\pi_{L}(m)}F_H \circ T_m \pi_{L}=T_m\pi _G$, from
expressions (\ref{tangent:1}) and (\ref{tangent:2}) we see that 
$T_{F_H([m])}\leaf=T_{F_H([m])} M_{\mathcal{O}_\rho}$, and the lemma holds.
\end{proof}

\begin{proof}[Proof of Theorem \ref{thm:symplectic:leaves}]
Note that
$M_{\mathcal{O}_\rho}=\mathcal{J}^{-1}(\mathcal{O}_{\rho})/G$ is
connected, so by the previous lemma it is an open subset of the
symplectic leaf $\leaf_{[m]}$. To see that
$M_{\mathcal{O}_\rho}=\leaf_{[m]}$ we note that $M
_{\mathcal{O}_\rho}$ is the accessible set of the projected flows of
the Hamiltonian vector fields associated to $G$-invariant
functions. But by (\ref{tangent:2}) so is $\leaf_{[m]}$ and hence the
equality holds.

The definition of $\omega_\rho$ shows that for any $f,h\in
C^{\infty}(M_H/L(H))$, we have
\begin{align*}
\{f,h\}_{M_H/L(H)}([m))&=\{f\circ\pi_L,h\circ\pi_L\}_{M_H}(m)
=\{\widetilde{f\circ\pi_L},\widetilde{h\circ\pi_L}\}_{M}(m)\\
&=\pi_\rho^*\omega_\rho(m)(X_{\widetilde{f\circ\pi_L}},
X_{\widetilde{h\circ\pi_L}})
=\omega_\rho([m])({X_f}|_{M_\rho},{X_h}|_{M_\rho}).
\end{align*}
Hence the symplectic leaves of $M_H/L(H)$ are the optimal
point reduced spaces $(\mathcal{J}^{-1}(\rho)/G_\rho,\omega_\rho)$. The
isomorphism $F_H:M_H/L(H)\to M_{(H)}/G$ now shows that the symplectic
leaves of $M_{(H)}/G$ are the optimal orbit reduced spaces $\left(\mathcal{J}^{-1}(\mathcal{O}_{\rho})/G,\omega_{\mathcal{O}_\rho}\right)$.
\end{proof}

\begin{remark}
Assume that the the original action $G\times M\to M$ 
is Hamiltonian with equivariant momentum map $\mu:M\to\gg^*$. Then the results
above yield the following:
\begin{enumerate}[(i)]
\item $M/G$ is a Poisson stratified space by orbit types $M_{(H)}/G$;
\item The reduced spaces $\mu^{-1}(\xi)/G_{\xi}$ are Poisson
stratified subspaces of $M/G$ (by orbit types).
\end{enumerate}
The singular spaces $\mu^{-1}(\xi)/G_{\xi}$ are not quotients
of smooth manifolds. Note that when the Poisson structure happens 
to be symplectic, the reduced spaces $\mu^{-1}(\xi)/G_{\xi}$ are 
symplectic stratified subspaces, but $M/G$ remains a Poisson stratified 
space: the strata of $\mu^{-1}(\xi)/G_{\xi}$
are the symplectic leaves of the strata of $M/G$.
\end{remark}

%%%%%%%%%%%%%%%%%%%%%%%%%%%%%%%%%%%%%%%%%%%%%%%%%
\subsection{An example}                         %
\label{subsec:example:Poisson:strat}                 %
%%%%%%%%%%%%%%%%%%%%%%%%%%%%%%%%%%%%%%%%%%%%%%%%%

Let $\Cc^{n+1}$ be the complex $(n+1)$-dimensional space with
holomorphic coordinates $(z_0,\dots,z_n)$ and anti-holomorphic
coordinates $(\overline{z}_0,\dots,\overline{z}_n)$. On the (real)
manifold $\Cc^{n+1}\setminus \{0\}$ we define a (real) quadratic Poisson
bracket by 
\[ \{z_i,z_j\}=a_{ij}z_iz_j,\] 
where $A=(a_{ij})$ is a fixed skew-symmetric matrix. The group $\Cc^*$
of non-zero complex numbers acts on $\Cc^{n+1}\setminus \{0\}$ by multiplication of
complex numbers. This is a free and proper Poisson action, so the
quotient $\CP(n)=\left(\Cc^{n+1}\setminus \{0\} \right)/\Cc^*$ inherits a Poisson bracket.

Let us consider now the action of the $n$-torus $\Tt^n$ on $\CP(n)$
defined by
\[ (\theta_1,\dots,\theta_n)\cdot [z_0:z_1:\cdots:z_n]=
[z_0,e^{i\theta_1}z_1,\cdots,e^{i\theta_n}z_n].\]
This is a Poisson action which is proper but not free. The quotient
$\CP(n)/\Tt^n$ can be identified with the standard simplex 
\[ 
\Delta^n=\left\{(\mu_0,\dots,\mu_n)\in{\mathbb R}^{n+1}\,\Bigg{|}\,\sum_{i=0}^n
\mu_i=1,\mu_i\ge 0\right\}.
\]
This identification is obtained via the map $\mu:\CP(n)\to\Delta^n$ defined by
\[ 
\mu([z_0:\cdots:z_n])=
\left(\frac{|z_0|^2}{|z_0|^2+\cdots+|z_n|^2},\cdots,\frac{|z_n|^2}{|z_0|^2+\cdots+|z_n|^2}\right).
\]
The strata of $\Delta^n=\CP(n)/\Tt^n$ are simply the faces of the
simplex of every dimension $0\le d\le n$.

Let us describe the Poisson nature of this stratification of
$\Delta^n=\CP(n)/\Tt^n$. The Poisson bracket on $\Delta^n$ is obtained
through the identification
\[ C^\infty(\Delta^n)\simeq C^\infty(\CP(n))^{\Tt^n}.\]
To see what it is, we simply determine the Poisson bracket between the
components of the map $\mu$. A straightforward computation yields
\begin{equation}
\label{eq:PoissonBr:simplex}
\{\mu_i,\mu_j\}_{\Delta}=\left(a_{ij}-\sum_{l=0}^n(a_{il}+a_{lj})\mu_l\right)\mu_i\mu_j,
\qquad (i,j=0,\dots,n).
\end{equation}
Now notice that (\ref{eq:PoissonBr:simplex}) actually defines a Poisson
bracket on ${\mathbb R}^{n+1}$. For this Poisson bracket, the interior of the 
simplex and its faces are Poisson submanifolds. A face 
$\Delta_{i_1,\dots,i_{n-d}}$ of dimension $0\le d\le n$ is given by
equations of the form:
\[ \sum_{i=0}^n \mu_i=1,\quad \mu_{i_1}=\cdots=\mu_{i_{n-d}}=0,
\quad\mu_i>0 \text{ for }i\not\in\{i_1,\dots,i_{n-d}\}.
\]
These equations define Poisson submanifolds since
\begin{enumerate}[(a)]
\item the bracket $\{\mu_i,\mu_l\}_{\Delta}$ vanishes whenever
$\mu_l=0$, and
\item the bracket $\{\mu_i,\sum_{l=0}^n\mu_l\}_{\Delta}$ vanishes
whenever $\sum_{l=0}^n\mu_l=1$.
\end{enumerate}
Therefore, the stratification of the simplex $\Delta^n$ by the faces 
is indeed a Poisson stratification.
\medskip

This example can be generalized in several directions. One can
consider, for example, more general homogeneous quadratic brackets
which are not necessarily holomorphic. Or one can consider other toric
manifolds, using Delzant's construction, which yield Poisson
stratifications of their Delzant polytopes.

%%%%%%%%%%%%%%%%%%%%%%%%%%%%%%%%%%%%%%%%%%%%%%%%
%%%%%%%%%%%%%%%%%%%%%%%%%%%%%%%%%%%%%%%%%%%%%%%%
\section{Momentum maps of Poisson actions}     %
\label{sec:momentum:maps}                      %
%%%%%%%%%%%%%%%%%%%%%%%%%%%%%%%%%%%%%%%%%%%%%%%%
%%%%%%%%%%%%%%%%%%%%%%%%%%%%%%%%%%%%%%%%%%%%%%%%

One can associate a canonical symplectic object to every 
Poisson manifold that can be thought of as its
symplectization. We show in this section that every Poisson action on a
Poisson manifold lifts to a globally Hamiltonian action on its
symplectization.

%%%%%%%%%%%%%%%%%%%%%%%%%%%%%%%%%%%%%%%%%%%%%%%%%%%%%%%%%
\subsection{Symplectization of a Poisson manifold}      %
\label{subsec:symplectization}                          %
%%%%%%%%%%%%%%%%%%%%%%%%%%%%%%%%%%%%%%%%%%%%%%%%%%%%%%%%%

Let $(M,\{\cdot,\cdot\})$ be a Poisson manifold with associated Poisson
tensor $\Pi$. We will denote by:
\begin{itemize}
\item $\X(M,\Pi):=\{X \in \X(M)\mid \Lie_X \Pi=0\}$ the Lie algebra of 
\emph{Poisson vector fields};
\item $\X_{\text{Ham}}(M,\Pi)\subset \X(M,\Pi)$ the Lie subalgebra 
of \emph{Hamiltonian vector fields}. 
\end{itemize}
There is a canonical symplectic object associated to the Poisson
manifold $(M,\Pi)$, namely, its \emph{Weinstein groupoid} $\Sigma(M)\tto M$
(\cite{CaFe,CrFe1,CrFe2}). We briefly  recall how this object is
defined.

A \textbf{cotangent path} in $M$ is a $C ^1 $ path $a:[0,1]\to T^*M$ such that
\begin{equation}
\label{eq:ctg:path}
\frac{\d}{\d t}p(a(t))=\sharp (a(t)),
\end{equation}
where $p:T^\ast M\to M$ is the canonical projection and
$\sharp:T^*M\to TM$ denotes the bundle map induced by the Poisson
tensor $\Pi$. The space of cotangent paths with the topology of uniform convergence will be denoted by
$P_\Pi(M)$. Notice that condition (\ref{eq:ctg:path}) defining a
cotangent path amounts to requiring the map $a\d t:TI\to T^*M$
to be a Lie algebroid morphism from the tangent Lie algebroid of the
interval $I:=[0,1]$ to the cotangent Lie algebroid $T^*M$ of the
Poisson manifold $M$. Given two cotangent paths $a_0,a_1\in P_\Pi(M)$
we say that they are \textbf{cotangent homotopic} if there exists a
family of cotangent paths $a_\eps\in P_\Pi(M)$ $(\eps\in [0,1])$,
joining $a_0$ to $a_1$, and satisfying the following property:
\begin{enumerate}
\item[(H)] For a connection $\nabla$ in $T^*M$ with torsion
  $T_{\nabla}$, the solution $b=b(\eps,t)$ of the differential
  equation
  \begin{equation*}
    %\label{eq:ctg:homotopy}
    \partial_{t}b-\partial_{\eps} a= T_{\nabla}(a, b),\qquad b(\eps,0)=0,
  \end{equation*}
  satisfies $b(\eps,1)=0$. 
\end{enumerate}
One can show that condition (H) is independent of the choice of
connection. This condition amounts to requiring the map 
\[ a\d t+b\d\epsilon:T(I\times I)\to T^*M \]
to be a Lie algebroid morphism. We
will write $a_0\sim a_1$ to denote that $a_0$ and $a_1$ are cotangent
homotopic paths. This is an equivalence relation on the set of
cotangent paths $P_\Pi(M)$.  For more details on cotangent paths and
homotopies we refer to \cite{CrFe1}.

The \textbf{Weinstein groupoid} $\Sigma(M)\tto M$ of the Poisson
manifold $(M,\{\cdot,\cdot\})$ is defined as follows:
\begin{enumerate}[(a)]
\item $\Sigma(M)$ is the space of equivalence classes of cotangent paths with the quotient topology:
  \[ \Sigma(M)=P_\Pi(M)/\sim ;\]
\item the source and target maps $\s, \t : \Sigma(M) \rightarrow  M $ are
  given by taking the initial and end-points of the paths:
  \[ \s([a])=p(a(0)),\quad \t([a])=p(a(1));\]
\item multiplication is given by concatenation of cotangent paths:
  \[ [a_1]\cdot [a_2]=[a_1\cdot a_2];\]
\item the identity section $i:M\to \Sigma(M)$ is obtained by taking the
  trivial cotangent path:
  \[ i(m)=[0_m],\quad (m\in M)\]
\item the inversion map $\iota:\Sigma(M)\to\Sigma(M)$ is obtained by 
  taking the opposite path:
  \[ \iota([a])=[\bar{a}],\]
  where $\bar{a}(t):= a(1-t)$.
\end{enumerate}
Note that $\Sigma(M)$ is a topological groupoid which is associated to
\emph{any} Poisson manifold. 
In the remainder of the paper we will work within the
class of \textbf{integrable Poisson manifolds} which means that
$\Sigma(M)$ is a Lie groupoid. A few exceptions to this assumption will be explicitly
noted.

The obstructions to integrability were determined in \cite{CrFe1,CrFe2}. When
$M$ is integrable, $\Sigma(M)\tto M$ is the unique source
simply-connected (i.e., the fibers of $\s$ are simply-connected)
Lie groupoid integrating the Lie algebroid $T^*M$. Moreover it is a
\emph{symplectic groupoid}: $\Sigma(M)$ carries a natural symplectic
2-form $\Omega\in\Omega^2(\Sigma(M))$ (\cite{CaFe,CrFe1}) which is
\emph{multiplicative}. We recall here the definition of multiplicative
form on a groupoid, since it will play an important role in the
sequel.

\begin{defn}
  Let $\G\tto M$ be a Lie groupoid. A differential form
  $\omega\in\Omega^\bullet(M)$ is called \textbf{multiplicative} if
  \[ m^*\omega=\pi_1^*\omega+\pi_2^*\omega,\]
  where $m:\G^{(2)}\to\G$ is the multiplication defined on the set of
  composable arrows $\G^{(2)}=\{(g,h) \mid \s(g)=\t(h)\}\subset\G\times\G$,
  and $\pi_1,\pi_2:\G^{(2)}\to\G$ are the (restrictions of the)
  projections onto the first and second factor, respectively.
\end{defn}

The multiplicative symplectic form $\Omega$ on $\Sigma(M)$ allows us to
identify the Lie algebroid $A=A(\Sigma(M))=\Ker T_M\s$ with the
cotangent Lie algebroid $T^*M$ via the isomorphism
\begin{equation}
\label{eq:iso:ctg:algebroid}
\Ker T_M \s\ni v\mapsto \left(i_v\Omega :TM\to\Rr\right).
\end{equation}
Here we identify $M$ with its image in $\Sigma(M)$ under the identity
section. For this and other basic properties, as well as a detailed
study of multiplicative 2-forms on Lie groupoids, we refer the reader to
\cite{BCWZ}.

In this paper, we would like to adopt the point of view that
$\Sigma(M)$ is \emph{the} symplectization of the Poisson manifold $M$.
For example, a basic fact is the following.

\begin{prop}
\label{prop:lifting:Poisson:maps}
Let $\phi:M\to M$ be a Poisson diffeomorphism. There exists a
symplectomorphism $\widetilde{\phi}:\Sigma(M)\to\Sigma(M)$ 
which covers $\phi$: it is the unique groupoid automorphism 
integrating the Lie algebroid automorphism $(T \phi^{-1})^*:T^*M
\rightarrow  T^*M$.
\end{prop}

The proof follows immediately from Lie's second theorem (which is
valid for Lie groupoids) since $\Sigma(M)$ is source
simply-connected. The explicit form of the map
$\widetilde{\phi}:\Sigma(M)\to\Sigma(M)$ is
\[ \widetilde{\phi}([a])=[(T\phi^{-1})^* \circ a], \]
for any cotangent path $a\in P_{\Pi}(M)$. In this form, we see that
$\widetilde{\phi}$ exists even in the non-integrable case. This is
one  instance that shows how advantageous it is having an explicit description
of $\Sigma(M)$ in terms of cotangent paths. We will see many other
examples later on.

A vector field $X\in\X(M)$ can be integrated over a cotangent path
$a\in P_\Pi(M)$ by setting
\[ \int_a X:=\int_0^1 \langle a(t),X(p(a(t)))\rangle \d t.\]
Note that for a Hamiltonian vector field $X_h\in\X_{\text{Ham}}(M,\Pi)$
the integral depends only on the end points
\[ \int_a X_h=h(p(a(1)))-h(p(a(0))).\]
A basic property which is proved in \cite{CrFe1} is the invariance of
the integral of Poisson vector fields under cotangent homotopies: if
$a_0,a_1\in P_\Pi(M)$ are cotangent homotopic paths and
$X\in\X(M,\Pi)$ is any Poisson vector field then
\[ \int_{a_0} X=\int_{a_1} X.\]
Therefore, we obtain a well defined map $c_X:\Sigma(M)\to\Rr$ by
setting
\[ c_X([a]):=\int_a X.\]
The additivity of the integral shows that $c_X$ is a groupoid 1-cocycle
\[ c_X([a_1]\cdot[a_2])=c_X([a_1])+c_X([a_2]);\]
equivalently, $c_X$ is a multiplicative 0-form.

On the other hand, the Van Est map associates to a
groupoid 1-cocycle $c:\Sigma(M)\to\Rr$ the Lie algebroid 1-cocycle
$\omega\in\Omega^1(A):=\Gamma(A^*)$ defined by
\[ \omega: m\mapsto T_m c|_{\Ker T_x s}, \quad \text{where} \quad x= 1_m.\]
Composing with the isomorphism $A\simeq T^*M$ given by 
(\ref{eq:iso:ctg:algebroid}) we obtain a Poisson vector field $X_c\in
\X(M,\Pi)$ (recall that Poisson vector fields are just Poisson cohomology cocycles or, which is the same, Lie algebroid cocycles for $T^*M$).
Summarizing, we have:
\begin{enumerate}[(i)]
\item the integration map, which associates to a Poisson vector field
  $X\in\X(M,\Pi)$ a (differentiable) groupoid cocycle $c_X\in
  C^1(\Sigma(M))$;
\item the Van Est map, which associates to a groupoid 1-cocycle $c\in
  C^1(\Sigma(M))$ a Poisson vector field $X_c\in \X(M,\Pi)$.
\end{enumerate}
The Van Est Theorem (see \cite{Cr1}) states that the correspondences $X\mapsto c_X$
and $c\mapsto X_c$ are inverses of each other. For details on these facts we refer 
the reader to \cite{CrFe1}.

%%%%%%%%%%%%%%%%%%%%%%%%%%%%%%%%%%%%%%%%%%%%%%%%%%%%%%%%%
\subsection{From Poisson actions to Hamiltonian actions}%
\label{subsec:lifting:action}                           %
%%%%%%%%%%%%%%%%%%%%%%%%%%%%%%%%%%%%%%%%%%%%%%%%%%%%%%%%%

Any Poisson action on $M$ can be lifted to an  action 
on $\Sigma(M)$ that admits a natural equivariant momentum map. 

\begin{thm}[Symplectization of Poisson actions]
\label{thm:lifted:action}
Let $G\times M\to M$ be a smooth action of a Lie group $G$ on a
Poisson manifold $M$ by Poisson diffeomorphisms. There exists a unique
lifted action of $G$ on $\Sigma(M)\tto M$ by symplectic groupoid
automorphisms. This lifted $G$-action is Hamiltonian and admits the
momentum map $J:\Sigma(M)\to\gg^*$ given by
\begin{equation}
\label{eq:lifted:momentum:map}
\langle J([a]),\xi\rangle=\int_a X_\xi,
\end{equation}
where $X_\xi\in\X(M,\Pi)$ denotes the infinitesimal generator determined by 
$\xi\in\gg$. Furthermore:
\begin{enumerate}[(i)]
\item The momentum map $J$ is $G$-equivariant and is a groupoid 1-cocycle.
\item The $G$-action on $M$ is Hamiltonian with momentum map $\mu:M\to\gg^*$
if and only if $J$ is an exact cocycle:
\[ J=\mu\circ\s-\mu\circ\t.\]
\end{enumerate}
\end{thm}

\begin{proof}
Applying Proposition \ref{prop:lifting:Poisson:maps} to each Poisson
automorphisms $\Phi_g:M\to M$, $g\in G$, defined by the Poisson action
$\Phi:G\times M\to M$, we obtain immediately a lifted symplectic
action $G\times \Sigma(M)\to\Sigma(M)$, with symplectic groupoid
automorphisms $\widetilde{\Phi}_g:\Sigma(M)\to\Sigma(M)$ that cover
$\Phi_g$:
\[
\xymatrix{
\Sigma(M)\ar[r]^{\widetilde{\Phi}_g}\ar@<.5ex>[d]\ar@<-.5ex>[d]&
\Sigma(M)\ar@<.5ex>[d]\ar@<-.5ex>[d]\\
M\ar[r]_{\Phi_g}&M}
\]

All we need to show is that the lifted action is Hamiltonian with
momentum map given by (\ref{eq:lifted:momentum:map}). Then the
remaining statements follow immediately from the expression of $J$. 

For each  $\xi\in\gg$, let $\widetilde{X}_\xi$ be the infinitesimal generator of the lifted action. It is a symplectic vector field and we need to
show that it is Hamiltonian, that is,
\begin{equation}
\label{eq:lifted:hamiltonian}
i_{\widetilde{X}_\xi}\Omega=\d J^\xi,
\end{equation}
where the Hamiltonian function $J^\xi$ is defined by
\[ J^\xi([a])=\langle J([a]),\xi\rangle.\]
We split the proof of (\ref{eq:lifted:hamiltonian}) into a few lemmas.

\begin{lem}
The 1-forms $i_{\widetilde{X}_\xi}\Omega$ and $\d J^\xi$ are multiplicative.
\end{lem}

\begin{proof}
First of all, $J^\xi$ is the groupoid cocycle that corresponds to the
Poisson vector field $X_\xi$. Hence, it is a multiplicative 0-form and
so its differential $\d J^\xi$ is a multiplicative 1-form. 

Now observe that the diagonal action of $G$ on $\Sigma(M)\times
\Sigma(M)$ has infinitesimal generator
$Y_{\xi}:=(\widetilde{X}_\xi,\widetilde{X}_\xi)$, leaves invariant the
space $\Sigma(M)^{(2)}$ of composable arrows, and makes the
projections $\pi_1,\pi_2:\Sigma(M)^{(2)}\to \Sigma(M)$
equivariant. Multiplication $m:\Sigma(M)^{(2)}\to \Sigma(M)$ is also
an equivariant map, since the action of $G$ on $\Sigma(M)$ is by
groupoid automorphisms. It follows that the infinitesimal generators
$Y_{\xi}$ and $\widetilde{X}_\xi$ are $\pi_1$, $\pi_2$, and
$m$-related. From this and using the fact that $\Omega$ is a
multiplicative 2-form, we see that
\begin{align*}
m^*\Omega&=\pi_1^*\Omega+\pi_2^*\Omega,\\
\Longrightarrow \quad
i_{Y_\xi}m^*\Omega&=i_{Y_\xi}(\pi_1^*\Omega+\pi_2^*\Omega),\\
\Longrightarrow \quad
m^*i_{\widetilde{X}_\xi}\Omega&=\pi_1^*i_{\widetilde{X}_\xi}\Omega
+\pi_2^*i_{\widetilde{X}_\xi}\Omega,
\end{align*}
so $i_{\widetilde{X}_\xi}\Omega$ is a multiplicative 1-form.
\end{proof}

\begin{lem}
The 1-forms $i_{\widetilde{X}_\xi}\Omega$ and $\d J^\xi$ are closed
and they coincide on $M$.
\end{lem}

\begin{proof}
$\d J^\xi$ is obviously closed. Also, since $\Omega$ is closed, we
find
\[ \d i_{\widetilde{X}_\xi}\Omega=\Lie_{\widetilde{X}_\xi}\Omega=0 .\]
To check that $i_{\widetilde{X}_\xi}\Omega$ and $\d J^\xi$ agree on $M$, we use the identification
\[ T_m\Sigma(M)= T_mM\oplus A_m\simeq T_mM\oplus T_m^*M,\]
provided by the isomorphism (\ref{eq:iso:ctg:algebroid}).
\end{proof}

Relation (\ref{eq:lifted:hamiltonian}) follows from the
previous two lemmas and the following result.

\begin{lem}
If two multiplicative 1-forms on a Lie groupoid $\G\tto M$ have the same
differential and agree on the identity section $M$ then they
must coincide.
\end{lem}

The proof of this lemma is exactly the same as in the case of
multiplicative 2-forms which is given in \cite[Corollary 3.4]{BCWZ},
so we omit it. This completes the proof of Theorem \ref{thm:lifted:action}.

\end{proof}
\vskip 15 pt

\begin{rems} $\quad$
  \begin{enumerate}[(i)]
  \item Since the action of $G$ on $\Sigma(M)$ is by groupoid
    automorphisms, all structure maps, i.e., the source and target maps
    $\s,\t:\Sigma(M)\to M$, the inversion map $\iota:\Sigma(M)\to\Sigma(M)$, and
    the identity section $i:M\to\Sigma(M)$, are $G$-equivariant.

  \item Theorem \ref{thm:lifted:action} is proved in \cite{CoDaWe,MiWe}
    for the special case of symplectic actions on symplectic manifolds (we
    will recover this case in the next section). In \cite{WeXu}, it is
    proved that every Poisson action lifts to a Hamiltonian action, but
    the explicit form of the momentum map is missing, since the description
    of the symplectic groupoid in terms of cotangent paths was not
    available.

  \item In \cite{MiWe}, the  authors consider group actions on symplectic
    groupoids by groupoid automorphisms. They show that if the groupoid is
    source simply-connected any such action has an equivariant momentum map
    which is a groupoid 1-cocycle. This follows also from Theorem
    \ref{thm:lifted:action} since any such action is the lift of a Poisson
    action and such a groupoid is isomorphic to the Weinstein groupoid
    $\Sigma(M)$.
  \end{enumerate}
\end{rems}

Note that if the original $G$-action on $M$ is Hamiltonian, so that 
$J:\Sigma(M)\to\gg^*$ is an exact 1-cocycle, then $J$ must 
vanish on the isotropy groups
\[ \Sigma(M,m)=\s^{-1}(m)\cap\t^{-1}(m).\]
In general, this not true and the restrictions
$J:\Sigma(M,m)\to\gg^*$
are non-trivial group homomorphisms.

\begin{defn}
The \textbf{group of periods} of the Poisson action at a point
$m\in M$ is the subgroup
\[ 
H_m=J(\Sigma(M,m))=\left\{\left. \xi\mapsto \int_a X_\xi \, \right| [a]\in\Sigma(M,m)\right\} 
\subset \mathfrak{g}^\ast.\]
\end{defn}

Therefore, the groups of periods of the action give natural obstructions
for a Poisson action to be a Hamiltonian action. We will see below
that, in the symplectic case, they are the only obstruction.

For distinct points $m_1,m_2\in M$, the groups of periods $H_{m_1}$ and 
$H_{m_2}$ are also distinct, in general. However, we have the following result.

\begin{prop}
If $m_1,m_2\in M$ are points that belong to the same symplectic leaf then their
groups of periods coincide, that is,  $H_{m_1}=H_{m_2}$.
\end{prop}

\begin{proof}
If $m_1,m_2\in M$ lie in the same symplectic leaf, we can find a cotangent
path $c:[0,1]\to M$ such that $\s([c])=m_1$ and $\t([c])=m_2$. Now, if $a$ is
a cotangent loop based at $m_1$, the concatenation $c\cdot a\cdot \bar{c}$ is
a cotangent loop based at $m_2$; here $\bar{c}$ is the oppositely oriented path $c $. The cocycle property of the momentum 
map gives
\[ J([c\cdot a\cdot \bar{c}])=J([c])+J([a])+J([\bar{c}])=J([a]).\]
This shows that $H_{m_1}\subset H_{m_2}$. Similarly, we have $H_{m_2}\subset H_{m_1}$, so
the result follows.
\end{proof}

%%%%%%%%%%%%%%%%%%%%%%%%%%%%%%%%%%%%
\subsection{Examples}              %
\label{subsec:examples}            %
%%%%%%%%%%%%%%%%%%%%%%%%%%%%%%%%%%%%

In this paragraph we illustrate Theorem \ref{thm:lifted:action} and
some its consequences by considering a few examples.

%%%%%%%%%%%%%%%%%%%%%%%%%%%%%%%%%%%%
\subsubsection{Symplectic actions} %
\label{subsec:symp:actions}        %
%%%%%%%%%%%%%%%%%%%%%%%%%%%%%%%%%%%%
Given a connected symplectic manifold $(M,\omega)$, so that $\Pi=\omega^{-1}$, 
the set of cotangent paths $P_\Pi(M)$ is naturally identified with the 
space of paths in $P(M)$: to a cotangent path $a:[0,1]\to T^*M$ we 
associate its base path $\gamma=p\circ a:[0,1]\to M$ and to a path 
$\gamma:[0,1]\to M$ we associate the cotangent path 
$a=(\sharp)^{-1}\dot{\gamma}$. Under this identification, a cotangent homotopy
becomes a standard homotopy (with fixed end points), and we conclude that
the Weinstein groupoid is the fundamental groupoid $\Sigma(M)=\Pi(M) $ of $M$.
The symplectic form $\Omega$ on $\Sigma(M)$ is given by
\[ \Omega=\s^*\omega-\t^*\omega,\]
and one checks immediately that it is multiplicative.

Applying Theorem \ref{thm:lifted:action} we recover the
following result of \cite{CoDaWe,MiWe}.

\begin{prop}
Let $G$ be a Lie group acting by symplectomorphisms on $(M,\omega)$. There
exists a unique lifted action of $G$ on $\Pi(M)$ by symplectic
groupoid automorphisms that covers the given action. The lifted
action is Hamiltonian with a momentum map $J:\Sigma(M)\to\gg^*$ given by
\[
\langle J([\gamma]),\xi\rangle=\int_\gamma i_{X_\xi}\omega,\quad \xi\in\gg.
\]
This map is a $G$-equivariant groupoid 1-cocycle.
\end{prop}

As we observed above, the original symplectic action is a Hamiltonian action
with momentum map $\mu:M\to\gg^*$ if and only if $J:\Pi(M)\to\gg^*$ is
an exact cocycle, in which case we have
\[ J=\mu\circ\s-\mu\circ\t.\]
In general, the symplectic action will not be Hamiltonian and this will
be reflected in the fact that the groups of periods are not trivial. These
are now given by
\[ 
H_{m}=J(\pi_1(M,m))=
\left\{\left. \xi\mapsto \int_\gamma i_{X_\xi}\omega\, \right| \, \gamma\in\pi_1(M,m)\right\} .
\]
In this case, there is only one symplectic leaf and the groups
of periods are all equal (the integral above only depends on 
the homology class of $\gamma$). Let us denote by $H\subset\gg^*$ 
this common group of periods. The composition
\[ \xymatrix{ \Pi(M)\ar[r]^{J}&\gg^*\ar[r]&\gg^*/H},\]
is now an exact groupoid 1-cocycle with values in the abelian group 
$\gg^*/H$ and it is the coboundary of a map 
$\mu:M \rightarrow \mathfrak{g}^\ast/H$. If we want $\mu:M\to \mathfrak{g}^\ast/H$ to 
be smooth, then we must replace $H$ by its closure $\overline{H}$. But then, since  $\overline{H}$ is a closed subgroup of $\gg^*$, the 
quotient $\gg^*/\overline{H}$ is isomorphic to a group $\Rr^s\times \Tt^t$, 
for some $s,t \in \mathbb{N}$, and we obtain
the \emph{moment r\'eduit} or \emph{cylinder valued momentum map} introduced by Condevaux-Dazord-Molino (\cite{CoDaMo}).

\begin{cor}
Let $G$ act by symplectomorphisms on $(M,\omega)$. Then the action admits
a momentum map $\mu:M\to\gg^*/\overline{H}$, in the sense that its fibers 
are preserved by the flows of the $G$-equivariant Hamiltonian vector fields.
\end{cor}

Note that our approach is more direct and even more canonical than the original 
one in \cite{CoDaMo} since it does not involve 
any choice of connection.

%%%%%%%%%%%%%%%%%%%%%%%%%%%%%%%%%%%%
\subsubsection{Cotangent bundles}  %
\label{subsec:ctg:bundle}          %
%%%%%%%%%%%%%%%%%%%%%%%%%%%%%%%%%%%%

In the previous paragraph we considered the symplectic case. At the
other extreme, we can consider manifolds with the zero Poisson
bracket, so that any action $G\times M\to M$ is a Poisson action. In
this case, the Weinstein groupoid of $M$ is the cotangent bundle
$\Sigma(M)=T^*M$, with its canonical symplectic form and multiplication
given by addition on the fibers; the source and target maps both
coincide with the canonical projection $T^*M\to M$. 

A cotangent path in $T^*M$ is just a path $a:[0,1]\to T_x^*M$ into a
fiber. Any such path is cotangent homotopic to a constant
$a(t)=\al_m\in T_m^*M$, namely its average
\[ \al_m=\int_0^1 a(t)\d t.\]
Therefore, it follows from the general expression (\ref{eq:lifted:momentum:map}) 
that the momentum map $j:T^*M\to\gg^*$ is given by
\begin{equation}
\label{eq:J:ctg:bundle}
\langle j(\al_m),\xi\rangle=\langle \al_m,X_\xi(m)\rangle,\quad
\al\in T_m^*M,\; \xi\in\gg.
\end{equation}

Therefore, in this case, Theorem \ref{thm:lifted:action} reduces to
the following well-known result.

\begin{prop}
Let $G$ be a Lie group acting on a manifold $M$. The cotangent lift of
$G$ to $T^*M$ is a Hamiltonian action with an equivariant momentum map
$j:T^*M\to\gg^*$ given by (\ref{eq:J:ctg:bundle}) and so it is linear on
the fibers.
\end{prop}

%%%%%%%%%%%%%%%%%%%%%%%%%%%%%%%%%%%%%%%%%%
\subsubsection{Linear and affine Poisson structures}%
\label{subsec:coadj:action}              %
%%%%%%%%%%%%%%%%%%%%%%%%%%%%%%%%%%%%%%%%%%

Let $\gg$ be any finite dimensional Lie algebra. The dual space $\gg^*$ carries
a natural linear Poisson structure $\Pi^{\text{lin}}$ and, conversely,
every linear Poisson structure on a vector space is of this form. More 
generally, one can consider an \emph{affine Poisson structure}
\[ \Pi=c+\Pi^{\text{lin}},\]
where $c\in\wedge^2\gg^*$ is a constant bivector field. If we view $c$
as a skew-symmetric bilinear map $c:\gg\to\gg^*$, then $[\Pi,\Pi]=0$
if and only if $c$ is a Lie algebra 1-cocycle with values in the coadjoint representation, that is, the following identity holds
\[ \langle\d c(\xi,\eta),\zeta\rangle=
	c([\xi,\eta],\zeta)+c([\eta,\zeta],\xi)+c([\zeta,\xi],\eta)=0,\quad \text{for all} \quad \eta,\xi,\zeta\in\gg.\]

Let us describe the symplectic groupoid $\Sigma(\gg^*_c)$ associated with the
Poisson structure $c + \Pi^{\text{lin}}$ on $\mathfrak{g}^\ast$. Let $G$ be the 1-connected Lie group with Lie algebra $\mathfrak{g}$. 
Then, as a manifold, $\Sigma(\gg^*_c)=T^*G$. If we identify $T^*G=G\times\gg^*$ using (say) 
left translations, then the source and target maps are given by
\[ \s(g,\alpha)=\alpha,\quad \t(g,\alpha)=\Ad^*g\cdot\alpha +C(g), \]
where $\Ad^*$ denotes the left coadjoint action and $C:G\to\gg^*$ is the group 
1-cocycle (with values in the coadjoint representation) integrating $c$. 
The multiplication is given by
\[ (g,\alpha)\cdot(h,\beta)=(gh,\beta),\]
provided $\alpha =\Ad^*g\cdot\beta+C(g)$, and the symplectic form 
on $\Sigma(\gg^*_c)$ is given by
\[ \Omega=\omega_\text{can}+B_c,\]
where $\omega_\text{can}$ is the standard symplectic form on the 
cotangent bundle $T^*G$, and $B_c$ is a magnetic term defined by
\[ B_c(g)(T_eL_g\cdot\xi,T_eL_g\cdot\eta)=-c(\xi,\eta),\quad \text{for all} \quad \xi,\eta\in\gg.\]
One easily checks that $\Omega$ is multiplicative and induces the given
affine Poisson structure on the base.

The $C$-twisted coadjoint action of $G$ on $\gg^*$ given by
\[ g\cdot\alpha:= \Ad^* g\cdot\alpha +C(g), \quad \text{for all} \quad 
g \in G, \; \alpha\in \mathfrak{g}^\ast,\]
is a Poisson action. Let us see what Theorem \ref{thm:lifted:action} 
yields in this case. Under the left trivialization of the
cotangent bundle, the lifted action of $G$ on $\Sigma(\gg^*_c)=T^*G$ 
is given by
\[ g\cdot(h,\beta)=(ghg^{-1},g\cdot\beta).\]
In fact, this formula defines an action by symplectic groupoid automorphisms
which restricts to the original action on the identity section. By uniqueness
of the lifted action, this must be the lifted action.

The expression of the equivariant momentum map $J:T^*G\to\gg^*$ for the lifted action, follows
from Theorem \ref{thm:lifted:action} and the fact that the original
action is actually a Hamiltonian action with momentum map
$\mu:\gg^*\to\gg^*$ the identity map. Hence, $J$ is an exact cocycle, that is, 
\[ J(g,\alpha)=\mu(\s(g,\alpha))-\mu(\t(g,\alpha))=\alpha-\Ad^* g\cdot\alpha-C(g) \quad \text{for all} \quad 
g \in G, \; \alpha\in \mathfrak{g}^\ast.\]

\begin{remark}
In the previous example, the symplectic groupoid $\Sigma(\gg^*_c)$ 
coincides with the action groupoid for the $C$-twisted coadjoint action.
\end{remark}

%%%%%%%%%%%%%%%%%%%%%%%%%%%%%%%%%%%%%%%%%%%%%%%%%%%%%
\section{Symplectization and reduction}             %
\label{sec:symp:reduction}                          %
%%%%%%%%%%%%%%%%%%%%%%%%%%%%%%%%%%%%%%%%%%%%%%%%%%%%%

In this section we investigate if the two basic operations of 
symplectization and reduction, discussed in the previous
two sections, commute. We start by looking at some general properties
of actions on symplectic groupoids, then we consider the case of free 
actions, and finally we consider the singular case.

%%%%%%%%%%%%%%%%%%%%%%%%%%%%%%%%%%%%%%%%%%%%%%%
\subsection{G-actions on symplectic groupoids}%
\label{subsec:actions:symp:grpd}              %
%%%%%%%%%%%%%%%%%%%%%%%%%%%%%%%%%%%%%%%%%%%%%%%

Let $\G\tto M$ be a symplectic groupoid with symplectic form $\Omega$,
so that the base is a Poisson manifold $(M,\{\cdot,\cdot\})$. We 
consider a Lie group $G$ that acts on $\G$ by symplectic groupoid
automorphisms with a momentum map $J:\G\to\gg^*$ satisfying:
\begin{itemize}
\item[(i)] $J$ is equivariant: $J(g\cdot x)=\Ad^*g\cdot J(x)$, for all $g\in G$ and $x\in\G$;
\item[(ii)] $J$ is a groupoid 1-cocycle: $J(x\cdot y)=J(x)+J(y)$, for all $x,y\in\G$ that
can be composed.
\end{itemize}
Note that the $G$-action on $\G$ covers a $G$-action on $M$ and this action is 
Poisson. The two assumptions above imply that the 0-level set 
\[ J^{-1}(0)\subset\G \]  
is a $G$-invariant subgroupoid over $M$. Note that condition (ii) above 
can also be expressed as saying that $J:\G\to \mathfrak{g}^*$ is a Lie groupoid 
homomorphism, where $\gg^*$ is viewed as an additive group (which is a groupoid 
over a one point space).

\begin{lem}
\label{lem:reduction}
Let $A=A(\G):=\Ker T_M\s$. As in {\rm (\ref{eq:iso:ctg:algebroid})}, $A$ is isomorphic to the
cotangent Lie algebroid $T^*M$ via the map provided by the symplectic form $\Omega$. Under the 
isomorphism $A\simeq T^*M$, the momentum map $j:T^*M\to \gg^*$ for the cotangent lifted action, 
given by {\rm (\ref{eq:J:ctg:bundle})}, corresponds to the Lie algebroid homomorphism $J_*:A\to \gg^*$ 
induced by $J:\G\to \mathfrak{g}^*$.
\end{lem}

\begin{proof}
Let $w\in A_m=\Ker T_{1_m} \s$. We find
\[ 
\langle J_*(w),\xi\rangle=\langle T_{1_m}J\cdot w,\xi\rangle
=\Omega(w,X_\xi(m))
\]
where the last equality is the momentum map condition for $J$. 
On the other hand, the isomorphism $A\simeq T^*M$ provided by $\Omega$,
maps the vector $w\in A_m$ to the covector $\al_m\in T_m^*M$ given by
$\langle\al_m,v\rangle=\Omega(w,v)$, for all $v\in T_mM$.
Therefore,  $J_*:A\to\gg^*$ corresponds to $j:T^*M\to\gg^*$ given 
by (\ref{eq:J:ctg:bundle}) under the isomorphism $A\simeq T^*M$.
\end{proof}

\begin{rem}
If $G$ is a connected Lie group, then an action $G\times M\to M$ is 
Poisson if and only if every infinitesimal generator $X_\xi$
is a Poisson vector field. This, in turn, is equivalent to the 
condition that $j:T^*M\to\gg^*$ is a Lie algebroid homomorphism.
\end{rem}

In particular, we conclude the following result.

\begin{cor}
Assume that $0$ is a regular value of $J$. Then $J^{-1}(0)\subset\G$ 
is a Lie subgroupoid which integrates the Lie subalgebroid
$j^{-1}(0)\subset T^*M$, 
where $j:T^*M\to \gg^*$ is the momentum map for the cotangent lifted
action {\rm (\ref{eq:J:ctg:bundle})}.
\end{cor}

Obviously, if the $G$-action on $\G$ is free then $0$ is a regular value of
the momentum map. Actually, this is equivalent to requiring that the 
$G$-action on the base is free, in view of the following proposition.

\begin{prop}
\label{prop:proper}
Let $G$ be a Lie group acting smoothly on a Lie groupoid $\G\tto M$ by
groupoid automorphisms. Then:
\begin{enumerate}[(i)]
\item The $G$-action on $\G$ is free if and only if the $G$-action on
$M$ is free.
\item The $G$-action on $\G$ is proper if and only if the $G$-action
on $M$ is proper.
\end{enumerate}
\end{prop}

\begin{proof}
To prove (i), note that if the $G$-action on $\G$ is free then 
the $G$-action on $M$ is obviously free, since the identity 
section embeds $M$ equivariantly into $\G$. For the converse we 
observe that the source map $\s:\G\to M$ is 
$G$-equivariant, so that we have $G_x\subset G_{\s(x)}$ for every $x\in
\G$.  Hence, if the $G$-action on $M$ is free so is the $G$-action on
$\G$.

Similarly, (ii) is trivially true in one direction, since the identity
section embeds $M$ equivariantly into $\G$. Conversely, assume that the
$G$-action on $M$ is proper. We want to show that the map
\[ G\times\G\to\G\times\G,\;\; (g,x)\mapsto (x,gx),\]
is proper. To see this,  pick sequences $\{g_k\}\subset\G$ and $\{x_k\}\subset\G$ 
such that $x_k\to x$ and $g_kx_k\to y$; we need to show that $\{g_k\}$ has
a convergent subsequence.  Setting $m_k=\s(x_k)$, $m=\s(x)$, and
$n=\s(y)$, the equivariance of $\s$ shows that $m_k\to m$ and
$g_km_k\to n$ in $M$.  Since the action of $G$ on $M$ is proper, we
conclude that $\{g_k\}$ has a convergent subsequence, as needed.
\end{proof}

We apply the prior results to the Weinstein groupoid of a 
Poisson manifold $M$.

\begin{prop}
\label{prop:free:lifted:action}
Let $G\times M \rightarrow M$ be a Poisson action. If $J:\Sigma(M)\to\gg^*$
is the momentum map {\rm (\ref{eq:lifted:momentum:map})} for the lifted action, 
then $J^{-1}(0)\subset\Sigma(M)$ is a subgroupoid. If the action is free 
(respectively, proper), then the lifted action 
$G\times\Sigma(M)\rightarrow\Sigma(M)$ is free (respectively, proper),
and $J^{-1}(0)\subset \Sigma(M)$ is a Lie subgroupoid integrating the Lie 
subalgebroid $j^{-1}(0)\subset T^*M$.
\end{prop}

%%%%%%%%%%%%%%%%%%%%%%%%%%%%%%%%%%%%
\subsection{The regular case}      %
\label{subsec:free:reduction}      %
%%%%%%%%%%%%%%%%%%%%%%%%%%%%%%%%%%%%

Let us look now at free and proper actions. Let $\G\tto M$ be a 
symplectic groupoid and let $G$ be a Lie group acting on $\G$ 
by symplectic groupoid automorphisms with equivariant momentum map 
$J:\G\to\gg^*$ which is also a groupoid 1-cocycle. If the action is
proper and free then we can form the symplectic quotient~(see \cite{mwr})
$\G/\!/G=J^{-1}(0)/G$.

\begin{prop}
If the action of $G$ on $\G$ is proper and free, then the 
symplectic reduced space 
\[ \G/\!/G:=J^{-1}(0)/G,\] 
is a symplectic groupoid which integrates the quotient
Poisson structure on $M/G$.
\end{prop}

\begin{proof}
Since $J:\G\to\gg^*$ is an equivariant groupoid 1-cocycle, the 0-level
set $J^{-1}(0)\subset\G$ is a Lie subgroupoid over $M$ which is
$G$-invariant. It follows that the symplectic quotient
\[ \G/\!/G=J^{-1}(0)/G \] 
carries both a groupoid structure over $M/G$ and a symplectic form
$\Omega_\text{red}$. The multiplicative property of $\Omega_\text{red}$
follows easily from the multiplicative property of $\Omega$.
 
It remains to show that the induced Poisson structure on the base
$M/G$ coincides with the quotient Poisson structure or, equivalently,
that the Lie algebroid of $\G/\!/G$ is isomorphic through $\Omega_\text{red}$
to the cotangent Lie algebroid $T^*(M/G)$. To check this, notice that
Lemma \ref{lem:reduction} shows that we have the following 
commutative diagram of Lie algebroids
\begin{equation}
\label{eq:diag:reduction}
\newdir{ (}{{}*!/-5pt/@_{(}}
\xymatrix{
T^*M\ar[d]& j^{-1}(0)\ar@<2pt>@{ (->}[l]\ar[r]\ar[d]& T^*(M/G)\ar[d]\\
A(\G)& A(J^{-1}(0))\ar@<2pt>@{ (->}[l]\ar[r]&A(\G/\!/G) }
\end{equation}
where the vertical arrows are the isomorphisms induced by $\Omega$ and 
$\Omega_\text{red}$. Therefore, $A(\G/\!/G)$ is isomorphic to $T^*(M/G)$.
\end{proof}

We now apply this proposition to the momentum map $J:\Sigma(M)\to\gg^*$ 
associated with a proper and free Poisson action. From 
the results in Section \ref{sec:momentum:maps} and Proposition 
\ref{prop:free:lifted:action} we conclude the following.

\begin{cor}
\label{cor:reduct:free:case}
Let $G\times M\to M$ be a proper and free Poisson action on an integrable
Poisson manifold $(M,\{\cdot,\cdot\})$. Then $M/G$ is an integrable 
Poisson manifold and the reduced symplectic groupoid
\[ \Sigma(M)/\!/G:=J^{-1}(0)/G,\]
is a symplectic groupoid integrating $M/G$.
\end{cor}

Note that, in general, the symplectic groupoid $\Sigma(M)/\!/G$ is not source
simply-connected, so it \emph{does not} coincide with $\Sigma(M/G)$.

\begin{ex}
\label{ex:counter:examp}
Let us consider the symplectic 4-manifold $M=\Cc^2 \setminus \{0\}$ with symplectic form
\[ \omega=\frac{i}{2}\left(\d z\wedge\d\overline{z}+\d w\wedge\d\overline{w}\right),\]
where $(z,w)$ are global holomorphic coordinates. On $M$ we consider the 
Hamiltonian $\Ss^1$-action
\[ \theta\cdot(z,w):=(e^{i\theta}z,e^{-i\theta}w), \quad \text{for all} \quad \theta\in\Ss^1\]
which has momentum map $\mu:M\to\Rr$ given by
\[ \mu(z,w)=\frac{1}{2}\left(|w|^2-|z|^2\right).\]

Since $M$ is simply-connected, its symplectic groupoid is the pair groupoid
\[ \Sigma(M)=M\times\overline{M},\quad \Omega=\omega\oplus(-\omega),\]
and the lifted action is given by
\[ \theta\cdot(z_1,w_1,z_2,w_2):=(e^{i\theta}z_1,e^{-i\theta}w_1,e^{i\theta}z_2,e^{-i\theta}w_2), \quad \text{for all} \quad  \theta\in\Ss^1.\]
This is again a Hamiltonian action with momentum map $J:\Sigma(M)\to\Rr$ given by
\begin{align*}
J(z_1,w_1,z_2,w_2)&=(\mu\circ\s-\mu\circ\t)(z_1,w_1,z_2,w_2)\\
&=\frac{1}{2}\left(|w_1|^2-|z_1|^2+|z_2|^2-|w_2|^2\right).
\end{align*}
We claim that the quotient symplectic groupoid $\Sigma(M)/\!/\Ss^1=J^{-1}(0)/ \Ss^1$
does not have source connected fibers.

First of all, we need to identify the quotient Poisson manifold $M/\Ss^1$. 
For that, we consider the following generating set of $\Ss^1$-invariant 
polynomials:
\begin{align*}
\sigma_1&=zw+\overline{zw},\\
\sigma_2&=i(zw-\overline{zw}),\\
\sigma_3&=|z|^2-|w|^2,\\
\sigma_4&=|z|^2+|w|^2,
\end{align*}
which satisfy the relations:
\[ \sigma_4^2=\sigma_1^2+\sigma_2^2+\sigma_3^2,\quad \sigma_4\ge 0.\]
Hence, the map $(\sigma_1,\sigma_2,\sigma_3):\Cc^2\to\Rr^3$
induces an isomorphism $M/\Ss^1\simeq \Rr^3 \setminus \{0\}$. For the
quotient Poisson structure, we just need to compute the Poisson brackets
\begin{align*}
\{\sigma_1,\sigma_2\}&=(\sigma_1^2+\sigma_2^2+\sigma_3^2)^{1/2}\\
\{\sigma_1,\sigma_3\}&=\{\sigma_2,\sigma_3\}=0,
\end{align*}
and these define the quotient Poisson structure if we view 
$(\sigma_1,\sigma_2,\sigma_3)$ as global coordinates on $M/\Ss^1$.
Therefore, the symplectic leaves of $M/\Ss^1$ are the horizontal planes 
$\{(\sigma_1,\sigma_2,\sigma_3)\mid \sigma_3=c\}$. Of course,
these symplectic leaves are just the reduced symplectic spaces $\mu^{-1}(c)/\Ss^1$.

We now look at the source fibers. The source maps form a commuting diagram
\[
\xymatrix{
J^{-1}(0)\ar[r]^{\widehat{\pi}}\ar[d]_{\s}& \Sigma(M)/\!/\Ss^1\ar[d]^\s\\
M\ar[r]_{\pi}& M/\Ss^1
}
\]
where on the horizontal rows we have principal $\Ss^1$-bundles. Therefore,
restricting to the $\s$-fibers, we get for each $[(z,w)]\in M/\Ss^1$ a
principal $\Ss^1$-bundle
\begin{equation}
\label{eq:cont:ex}
\widehat{\pi}:\s^{-1}(\Ss^1\cdot(z,w))\to\s^{-1}([(z,w)]).
\end{equation}
For example, we have
\[ 
\s^{-1}(\Ss^1\cdot(z,0))=\{(|z|e^{i\theta},0,z_2,w_2)\mid |w_2|^2-|z_2|^2=|z|^2\}\simeq \Ss^1\times\Ss^1\times\Cc.
\]
Hence, for the point $[(z,0)]$, the first terms in the long exact homotopy sequence of (\ref{eq:cont:ex})
are:
\[
\xymatrix{\ar[r]&\pi_2(\s^{-1}([(z,0)])\ar[r]&\Zz\ar[r]&\Zz\times\Zz\ar[r]&\pi_1(\s^{-1}([(z,0)])\ar[r]&0}
\]
which shows that $\s^{-1}([(z,0)])$ is not simply connected.
\end{ex}

Our task now is to understand how symplectization and reduction
fail to commute. First of all, $J^{-1}(0)$ (and, hence, also its source fibers) 
need not be connected, so we consider its connected component containing the identity
section, denoted $J^{-1}(0)^0$. This is a subgroupoid of $\Sigma(M)$ which, in the integrable
case, has the same Lie algebroid as $J^{-1}(0)$. 

\begin{prop}
$J^{-1}(0)^0\subset \Sigma(M)$ is the unique source connected Lie subgroupoid
integrating the Lie algebroid $j^{-1}(0)\subset T^*M$. In particular,
\begin{equation} 
\label{eq:conn:subgrpd}
J^{-1}(0)^0=\{[a]\in\Sigma(M) \mid j(a(t))=0,\,\forall t\in[0,1]\}.
\end{equation}
\end{prop}

\begin{proof}
As we saw in Lemma \ref{lem:reduction}, the groupoid homomorphism $J:\Sigma(M)\to\gg^*$
corresponds to the Lie algebroid homomorphism $j:T^*M\to\gg^*$. This implies that
$J^{-1}(0)$, and hence also $J^{-1}(0)^0$, has Lie algebroid $j^{-1}(0)$. Now, both 
sides of (\ref{eq:conn:subgrpd}) are source connected
Lie subgroupoids of $\Sigma(M)$ integrating the subalgebroid $j^{-1}(0)\subset T^*M$, hence
they must coincide.
\end{proof}

If $G$ is connected,  $J^{-1}(0)^0$ is $G$-invariant. In the sequel we will 
consider only the case where $M$ is integrable and $G$ is connected. Consider 
then the groupoid diagram
\[
\newdir{ (}{{}*!/-5pt/@_{(}}
\xymatrix{
\Sigma(M)\ar@<3pt>[d]\ar@<-3pt>[d]& J^{-1}(0)^0\ar@<2pt>@{ (->}[l]\ar[r]^{\Phi}\ar@<3pt>[d]\ar@<-3pt>[d]& J^{-1}(0)^0/G\ar@<3pt>[d]\ar@<-3pt>[d]\\
M\ar@{=}[r]& M\ar[r]& M/G}
\]
with an associated Lie algebroid diagram (see Lemma \ref{lem:reduction} 
and the commutative diagram (\ref{eq:diag:reduction}))
\[
\newdir{ (}{{}*!/-5pt/@_{(}}
\xymatrix{
T^*M\ar[d]& j^{-1}(0)\ar@<2pt>@{ (->}[l]\ar[r]^{\phi}\ar[d]& T^*(M/G)\ar[d]\\
M\ar@{=}[r]& M\ar[r]& M/G.}
\]
Denote by $\G(j^{-1}(0)))$ the source 1-connected groupoid integrating $j^{-1}(0)$. 
The Lie algebroid morphism $\phi$ integrates to a morphism of source 1-connected groupoids
\[
\xymatrix{
\G(j^{-1}(0))\ar[r]^{\widehat{\Phi}}\ar@<3pt>[d]\ar@<-3pt>[d]& \Sigma(M/G)\ar@<3pt>[d]\ar@<-3pt>[d]\\
M\ar[r]& M/G}
\]
which covers the homomorphism $\Phi$
\begin{equation}
\label{eq:diag:com:symp}
\newdir{ (}{{}*!/-5pt/@^{(}}
\xymatrix{
K_M\ar[d]^{\widehat{\Phi}}\ar@<-2pt>@{ (->}[r]& \G(j^{-1}(0))\ar[d]^{\widehat{\Phi}}\ar[r]^{\widehat{p}}& J^{-1}(0)^0\ar[d]^{\Phi}\\
K_{M/G}\ar@<-2pt>@{ (->}[r]& \Sigma(M/G)\ar[r]^p& J^{-1}(0)^0/G.}
\end{equation}
Note that $K_M$ and $K_{M/G}$ are group bundles over $M$ and $M/G$, respectively,  whose fibers are discrete.

\begin{lem}
For every $m\in M$, the map $\widehat{\Phi}$ induces an 
isomorphism between $K_m$ and $K_{G\cdot m}$.
\end{lem}

\begin{proof}
In the commutative diagram (\ref{eq:diag:com:symp}) 
the maps $\Phi$ and $\widehat{\Phi}$ are just principal bundle projections for
the (proper and free) action of $G$. In fact, the action of 
$G$ on $\G(j^{-1}(0))$ comes from the identification of this
groupoid as equivalence classes of $A$-paths for the Lie algebroid
$A=j^{-1}(0)$. But an $A$-path $a:I\to j^{-1}(0)$ is just a cotangent
path with values in $j^{-1}(0)$, while an $A$-homotopy is 
just a cotangent homotopy with values in $j^{-1}(0)$. Since
$j^{-1}(0)$ is $G$-invariant, the $G$-action on cotangent paths gives 
an action of $A$-paths which preserves $A$-homotopies. This action
is proper and free, just like the lifted action on $\Sigma(M)$.

It follows that the action of $G$ on $K_M$ sends the fiber $K_m$
isomorphically to the fiber $K_{g\cdot m}$ for each $m\in M$ and $g\in G$.
Hence the the map $\widehat{\Phi}$ induces an isomorphism between 
$K_m$ and $K_{G\cdot m}$.
\end{proof}

Therefore, we see that the group bundle $K_M$ measures how
symplectization and reduction fail to commute. The fibers of this
bundle have a simple description in terms of cotangent paths.

\begin{thm}[Symplectization commutes with reduction]
  \label{thm:symp:reduct:alt}
  Let $G\times M\to M$ be a proper and free Poisson action. Then symplectization 
  and reduction commute if and only if the following groups
  \[ 
  K_{m}:=\frac{\{a:I\to j^{-1}(0) \mid a \text{\rm ~is a cotangent loop such that }a\sim 0_m\}}
   {\{\text{\rm cotangent homotopies with values in } j^{-1}(0)\}}
  \]
  are trivial, for all $m\in M$.
\end{thm}

\begin{proof}
We just need to combine the commutative diagram (\ref{eq:diag:com:symp}) with
the previous lemma and the description of $\G(j^{-1}(0))$ as cotangent paths to conclude that 
  \[ 
  \G(j^{-1}(0))=\frac{\{a:I\to j^{-1}(0) \mid a \text{ is a cotangent path}\}}
 {\{\text{cotangent homotopies with values in } j^{-1}(0)\}}
  \]
which immediately implies the statement in the theorem.
\end{proof}

Note that for the trivial Poisson bracket the cotangent paths lie in the fibers of the projection and the cotangent homotopies $T^*M\to M$ are just homotopies (with free end points). Hence, the
groups $K_m$ are all trivial and the theorem reduces to the following well-known result on cotangent bundle reduction~\cite{fom}.

\begin{cor}
\label{cor:cotg:reduct}
For a proper and free action $G\times M\to M$, we have
\[ T^*(M/G)=T^*M/\!/G. \]
\end{cor}

In the symplectic case, we can identify cotangent paths with ordinary paths (see
Section \ref{subsec:symp:actions}), so in the case of Hamiltonian $G$-spaces 
the groups $K_m$ take a special simple form and they can be described in terms 
of the fibers of the momentum map.

\begin{cor}
\label{cor:symp:reduct:alt}
Let $G\times M\to M$ be a Hamiltonian action on a symplectic manifold $(M,\omega)$
with momentum map $\mu:M\to\gg^*$. Then symplectization and reduction commute 
if and only if the following groups
\[ K_{m}:=\Ker i_*\subset \pi_1(\mu^{-1}(c),m)\]
are trivial for all $m\in M$, where $c=\mu(m)$ and $i:\mu^{-1}(c)\hookrightarrow M$ is the 
inclusion. 
\end{cor}

The homotopy long exact sequence of the pair $(M,\mu^{-1}(c))$ gives the   exact
sequence
\[
\xymatrix{ 
\pi_2(M,\mu^{-1}(c),m)\ar[r]^{\partial}&\pi_1(\mu^{-1}(c),m)\ar[r]^{i_*}&
\pi_1(M,m)\ar[r]^{j_*\qquad}&\pi_1(M,\mu^{-1}(c),m)}.
\]
Hence, for a Hamiltonian action, the groups vanish if the fibers of the 
momentum map are simply connected, or if the second relative homotopy groups 
of the fibers vanish. The latter occurs when the group is compact and the momentum map 
is proper as the following result shows.

\begin{cor}
\label{cor:symp:reduct:alt:2}
Let $G\times M\to M$ be a free Hamiltonian action of a compact Lie group on a 
symplectic manifold $(M,\omega)$ with a proper momentum map $\mu:M\to\gg^*$. Then
symplectization and reduction commute. Moreover, the isotropy groups $(\Sigma(M)/\!/G)_{[m]}$
all have the same number of connected components, that is, 
\[ 
\pi_0\left((\Sigma(M)/\!/G)_{[m]}\right) =\pi_1(M/G)=\pi_1\left(M_{\red},[m]\right), 
\]
where $M_{\red}=\mu^{-1}(\mathcal{O}_c)/G$ is the symplectic reduced space at value $c=\mu(m)$.
\end{cor}

\begin{proof}
Let us start by observing that it is enough to prove the result for the symplectic reduced space
at zero: $M/\!/G:=\mu^{-1}(0)/G$. Indeed, using the standard shifting trick, we can identify the symplectic
reduced space $M_{\red}:=\mu^{-1}(\mathcal{O}_c)/G$, at an arbitrary value $c\in\gg^*$, with $(M\times \mathcal{O}_c)/\!/G$,
the symplectic quotient at $0$ for the diagonal Hamiltonian $G$-action on $M\times\overline{\mathcal{O}_c}$, where $ \overline{\mathcal{O}_c} $  denotes the coadjoint orbit through $ c \in  \mathfrak{g}^\ast $  endowed with the negative orbit symplectic form. A momentum 
map for this action is $M\times\overline{\mathcal{O}_c}\ni(m,\xi)\mapsto \mu(m)-\xi\in\gg^*$. Since $\mathcal{O}_c$ is compact, 
this is a proper momentum map and, if the result holds for the symplectic reduced spaces at zero, we conclude that
\[ \pi_1(M_{\red})=\pi_1((M\times \mathcal{O}_c)/G). \]
Since coadjoint orbits of compact Lie groups are 1-connected we obtain $\pi_1(M_{\red})=\pi_1(M/G)$,
as claimed.

Now, since $\mu:M\to\gg^*$ is a proper submersion, by a classical result of Ehresmann, it is a locally trivial 
fibration with fiber type $\mu^{-1}(0)$. Moreover, the fibers are connected (see \cite{LMTW}) and the base is contractible. This shows that 
the map $i_*:\pi_*(\mu^{-1}(0))\to \pi_*(M)$ is an isomorphism. Now the inclusion of fibrations
\[
\newdir{ (}{{}*!/-5pt/@^{(}}
\xymatrix{
\mu^{-1}(0)\ar[d] ^{\quad}="a"&M\ar[d]_{\quad}="b"\\
\mu^{-1}(0)/G &M/G \ar @{ (->} "a";"b"
}
\]
leads to the commutative diagram:
\[
\xymatrix{
\pi_1(G)\ar[r]\ar[d] &\pi_1(\mu^{-1}(0))\ar[r]\ar[d]&\pi_1(\mu^{-1}(0)/G)\ar[r]\ar[d]&\pi_0(G)\ar[r]\ar[d] &\pi_0(\mu^{-1}(0))\ar[d]\\
\pi_1(G)\ar[r] &\pi_1(M)\ar[r]&\pi_1(M/G)\ar[r]&\pi_0(G)\ar[r] &\pi_0(M)
}
\]
The horizontal lines are exact and we already know that all vertical arrows except the middle one are isomorphisms.
By the Five Lemma, the middle one is also an isomorphism and we conclude that:
\[ \pi_1(M/\!/G)=\pi_1(\mu^{-1}(0)/G)=\pi_1(M/G). \]

Since $\Ker i_*$ is trivial, Corollary \ref{cor:symp:reduct:alt} shows that reduction and symplectization commute:
\[ \Sigma(M/G)=\Sigma(M)/\!/G.\]
Now note that $M/\!/G$ is a symplectic leaf of $M/G$, so it is an orbit of this groupoid. For any 
$[m]\in M/\!/G$, the homotopy long exact sequence of the target fibration $\t:\s^{-1}([m])\to M/\!/G$ gives
\[ 
\pi_1\left( M/\!/G,[m] \right) =\pi_0\left( \t^{-1}([m])\cap\s^{-1}([m]) \right) =\pi_0\left( \Sigma(M/G)_{[m]} \right) =\pi_0\left(\Sigma(M)/\!/G)_{[m]}\right), \]
which proves the corollary.
\end{proof}

Note that, in general, we \emph{do not} have $\pi_1(M)=\pi(M_{\red})$, contrary to what happens for 
Hamiltonian actions of compact Lie groups on compact symplectic manifolds (see \cite{Li}). A simple example
is provided by the cotangent lifted action of a compact, connected, non simply-connected Lie group $G$ on $M=T^*G$.
All assumptions of Corollary \ref{cor:symp:reduct:alt:2} are satisfied and we have 
$\pi_1(M_{\red})=\pi_1(M/G)=\pi_1(\gg^*)=1$, but $\pi_1(M)=\pi_1(T^*G)=\pi_1(G)\not=1$.

\begin{example}
Using Corollary \ref{cor:symp:reduct:alt} it is easy to check that ``symplectization
commutes with reduction'' fails in Example \ref{ex:counter:examp}. In this 
example, the fibers of the momentum map $\mu:\Cc^2 \setminus \{0\}\to\Rr$ have diffeomorphism type
\[
\mu^{-1}(c)\simeq
\begin{cases}
\Cc\times\Ss^1,\text{ if }c\ne 0,\\
(\Cc\setminus \{0\})\times\Ss^1,\text{ if }c=0.
\end{cases}
\]
Since $M$ is simply connected, we find that 
\[
K_m\simeq\pi_1(\mu^{-1}(c))=
\begin{cases}
\Zz,\text{ if }c\ne 0,\\
\Zz\times\Zz,\text{ if }c=0,
\end{cases}
\]
and so are non-trivial. 

Note also that in this example the momentum map is not proper and Corollary \ref{cor:symp:reduct:alt:2} does not hold. 
In fact, $\pi_1(M/G)$ is trivial, the symplectic reduced spaces at non-zero level also 
have trivial fundamental group, but at zero level we have $\pi_1(M/\!/G)=\Zz$.
\end{example}

%%%%%%%%%%%%%%%%%%%%%%%%%%%%%%%%%%%%
\subsection{The singular case}     %
\label{subsec:singular:reduction}  %
%%%%%%%%%%%%%%%%%%%%%%%%%%%%%%%%%%%%

Let us now drop the assumption that the action is free, so we let
$G\times M\to M$ be a proper action of a Lie group $G$ on a Poisson
manifold $M$ by Poisson diffeomorphisms. We continue to assume that
$(M,\{\cdot,\cdot\})$ is an integrable Poisson manifold. We will
consider the orbit type decomposition of $M$ and $X:=M/G$ and, as before, 
we will ignore connectedness issues by  working with $\Sigma$-manifolds, if necessary. 
Assuming that $G$ is connected, the strata of $M$ and $X$ are, respectively, 
$M_{(H)}$ and $M_{(H)}/G$.

First of all, by Proposition \ref{prop:proper}, the lifted $G$-action on 
$\Sigma(M)$ is proper and, by Theorem \ref{thm:lifted:action}, it is 
a Hamiltonian action with an equivariant momentum map $J:\Sigma(M)\to\gg^*$. We can
form the symplectic quotient $ J^{-1}(0)/G$
which is now a symplectic stratified space (\cite{LeSj, OrRa}). Note that:
\begin{itemize}
\item since $J$ is a groupoid 1-cocycle, $J^{-1}(0)\subset\Sigma(M)$ is still
a subgroupoid, but it \emph{fails} to be a Lie groupoid, since $0$ is not
a regular value;
\item since the action of $G$ on $J^{-1}(0)$ is not free, the groupoid 
multiplication \emph{does not} drop to $J^{-1}(0)/G$. 
\end{itemize}
We need to fix this unhappy state of affairs. 

Our main observation is the following. Instead of thinking of $M$ as a \emph{smooth} 
manifold, one should think of it as a \emph{stratified} space. This means that
one should replace $\Sigma(M)$ by a groupoid integrating the stratified
space $M=\bigcup_{(H)}M_{(H)}$. Therefore, one must deal with the following 
two issues:
\begin{enumerate}
\item[(i)] One must consider stratified Lie groupoids and Lie algebroids.
\item[(ii)] The strata $M_{(H)}$ are not Poisson manifolds, but rather Dirac manifolds,
so $M$ has a Dirac stratification.
\end{enumerate}
Once these two difficulties are solved, we will see that the results we have obtained
for the free case, extend to the non-free case.

\subsubsection{Stratified Lie theory}

We will not develop here a complete stratified Lie theory; we limit ourselves to 
outline some results leaving the details for future work. In what follows we will 
assume that all stratifications are {\bf smooth stratifications} (see, e.g., \cite{Pflaum}), that is, we assume the existence of a stratified atlas formed by singular charts that can then be used to define the set of $ C^\infty $ functions. A continuous map between two smooth stratified spaces is called a {\bf smooth map} if it pulls back smooth functions to smooth functions. If a smooth map between stratified spaces is, additionally, a morphism of stratified spaces, we say that it is a {\bf smooth stratified map} (note that a smooth map does not need to be a stratified map and, conversely, a stratified map is not necessarily smooth). For example, a smooth stratified map $\phi: (X _1, \mathcal{S} _1) \rightarrow (X _2, \mathcal{S} _2) $ is an {\bf immersion} (respectively, {\bf submersion}) if  $\phi $ is smooth and for any stratum $S  \in 
\mathcal{S} _1 $ there exists a stratum $T \in \mathcal{S} _2$ such that  $\phi (S)\subset T $ and  $\phi|_{S}:S \rightarrow  T $ is a smooth immersion (respectively, submersion).

Note that it maybe cumbersome to check that a given stratification is smooth, but in the
case we are mostly interested (orbit type stratifications associated with proper group actions) this is a straightforward consequence of the Slice Theorem which, additionally shows that the stratification is Whitney.

\begin{defn}
A \textbf{stratified Lie groupoid} is a groupoid in the category of Whitney smooth stratified spaces.
\end{defn}

Hence, both the space of arrows $\G$ and the space of objects $X$ are smooth stratified
spaces and all structure maps (source, target, identity, inversion, and multiplication) 
are smooth stratified maps. Moreover, $\s$ and $\t$ are submersions and the identity section is an embedding. The Whitney condition ensures that the tangent bundle of the stratified Lie groupoid is also a smooth stratified space (see~\cite{Pflaum} for more details).

Similarly, we have the notion of a stratified Lie algebroid over a Whitney smooth stratifed space.

\begin{defn}
A \textbf{stratified Lie algebroid} is a smooth stratified vector bundle $A\to X$, with a Lie bracket $[~,~]$ on the space of smooth sections, and a smooth stratified vector bundle map $\rho:A\to TX$ satisfying the Leibniz identity
\[ [\al,f\be]=f[\al,\be]+\rho(\al)(f)\be,\quad (\al,\be\in\Gamma(A), f\in C^\infty(X)), \]
where $TX$ denotes the stratified tangent bundle.
\end{defn}

We leave to the reader the task of defining morphisms of stratified Lie groupoids and 
stratified Lie algebroids. 

\begin{prop}
Every stratified Lie groupoid $\G\tto X$ has an associated stratified 
Lie algebroid $A\to X$.
\end{prop}

\begin{proof}
The proof is entirely analogous to the classical case. Just observe that $A:=\Ker d\s|_{X}$ is a
stratified vector bundle and then that $d\t$ restricts to a stratified vector bundle map 
from $A$ to $TX$. The Lie bracket on sections of $A$ is deduced by exactly the same procedure as in the smooth case.  Note that $A \subset T \G$ is by construction a stratified subspace.
\end{proof}

Given a stratified Lie algebroid one can construct the associated Weinstein groupoid
exactly by the same methods as in \cite{CrFe2}. In the special case of the stratified
cotangent Lie algebroid associated with a Poisson stratification we obtain a stratified
symplectic Lie groupoid.

\begin{example}
Let us return to the Poisson stratification of the simplex $\Delta^n$ constructed in Section \ref{subsec:example:Poisson:strat}. This defines a Lie algebroid structure on the stratified cotangent bundle $T^*\Delta^n$. The corresponding stratified Lie groupoid coincides with the restriction to $\Delta^n$ of the (smooth) Lie groupoid integrating $\Rr^{n+1}$ with the cubic Poisson structure (\ref{eq:PoissonBr:simplex}).
\end{example}

\subsubsection{The Dirac stratification}

Let now $G\times M\to M$ be a proper Poisson action on a Poisson manifold $(M,\Pi)$. 
As before, we consider the orbit type decomposition of $M$ and $X=M/G$. As before, we ignore
connectedness issues by working with $\Sigma$-manifolds, so that the strata are respectively, $M_{(H)}$ and $M_{(H)}/G$.

Consider the lifted cotangent bundle action $G\times T^*M\to T^*M$. As we have noticed in
Section \ref{subsec:Poisson:stratifications}, each orbit type manifold $M_{(H)}$ has
an induced Dirac structure $L_{(H)}\subset TM_{(H)}\oplus T^*M_{(H)}$. Now consider
the \textbf{stratified generalized tangent bundle}:
\[ \T M=\bigcup_{(H)} TM_{(H)}\oplus T^*M_{(H)}\]
(since the generalized tangent bundle is also known to people working in control theory as the \emph{Pontryagin bundle}, perhaps a better name to denote $\T M$ would be the \emph{stratified Pontryagin bundle}). We can put together the Dirac structures on each stratum, obtaining a stratified Dirac structure
\[ L=\bigcup_{(H)} L_{(H)}\subset \T M.\]
Note that $L\to M$ is a stratified Lie algebroid. If $(M,\Pi)$ is integrable
then each $L_{(H)}$ is integrable. 

\begin{rem} 
Here we are using the fact that the orbit type stratification 
of $M$ is Whitney. This property guarantees that $\T M $ is also a smooth stratified space. Then
$L \subset \T M $ is a smooth stratified subspace.
\end{rem}

\begin{prop}
Let $G\times M\to M$ be a proper Poisson action on an integrable Poisson manifold $(M,\Pi)$. 
The Weinstein groupoid $\Sigma(M,L)$ of the stratified Dirac structure $L$ associated with $\Pi$
is a stratified presymplectic groupoid integrating $L$. 
\end{prop}

\begin{proof}[Sketch of proof]
The groupoid $\Sigma(M,L)$ is defined by the same method used in the smooth case (see \cite{CrFe2}). Note that the $A$-paths take values in each $L_{(H)}$ (i.e., they don't cross
different strata of $L$), so that over each $L_{(H)}$ we obtain a smooth presymplectic groupoid as in \cite{BCWZ} (recall that, by assumption, each $L_{(H)}$ is integrable). Hence, $\Sigma(M,L)$
becomes a smooth stratified presymplectic Lie groupoid. According to our definition of stratified Lie groupoid one still needs to check that the Whitney condition is satisfied.
\end{proof}

\subsubsection{Symplectization and reduction in the non-free case}
The action $G\times M\to M$ preserves the strata $M_{(H)}$, so there is a lifted action on 
the stratified generalized tangent bundle $\T M$ (which assembles together the tangent and cotangent lifted actions). This lifted action preserves the Dirac structures $L_{(H)}$, so it lifts to an action on the stratified presymplectic groupoid $\Sigma(M,L)$ by symplectic groupoid
automorphisms.

On the other hand, the stratified map $j:L\to\gg^*$ given by
\[ j(v,\al)=j_{(H)}(\al),\]
where $j_{(H)}:T^*M_{(H)}\to\gg^*$ is the momentum map for the cotangent lifted action
of $G$ on $T^*M_{(H)}$, is a morphism of Lie algebroids. Hence, it lifts to a stratified
map $J:\Sigma(M,L)\to\gg^*$ which is a morphism of Lie groupoids. One checks easily that
$J^{-1}(0)\subset \Sigma(M,L)$ is a stratified Lie subgroupoid integrating the stratified
Lie subalgebroid $j^{-1}(0)\subset L$ and we can introduce the quotient.
Now, just like in the free case, we have the following result.

\begin{thm}
Let $G\times M\to M$ be a proper action on a integrable Poisson manifold. Then $J^{-1}(0)\subset \Sigma(M,L)$ is a 
stratified Lie subgroupoid with stratified Lie algebroid $j^{-1}(0)\subset L$.
The quotient
\[ \Sigma(M,L)/\!/G:=J^{-1}(0)/G \]
is a stratified symplectic groupoid integrating the Poisson stratification $X=M/G$.
\end{thm}

\begin{proof}
The first statement was already proved. It is also clear that $\Sigma(M,L)/\!/G$ is a stratified Lie groupoid. 

In order to check that there is a reduced symplectic form $\Omega_{\text{red}}$, we can restrict to a single stratum $M_{(H)}$ where there is a single orbit type. Then $M_{(H)}\to M_{(H)}/G$ is a smooth bundle,
and all results for free actions apply. In particular, we have the commutative diagram
of Dirac structures (see the discussion before Theorem \ref{Poisson strata we go}), which 
makes it clear that the corresponding diagram of Lie algebroids commutes. This implies that $\Sigma(M,L)/\!/G$ is symplectic.

Since $G$ acts by groupoid automorphisms and the presymplectic form $\Omega$ is multiplicative,
it is clear that $\Omega_{\text{red}}$ is also multiplicative, and hence $\Sigma(M,L)/\!/G$ becomes a stratified symplectic groupoid.
\end{proof}

Note that the stratified Lie algebroid $\Sigma(M,L)/\!/G$ is the stratified cotangent bundle $T^*X$
with the Lie algebroid structure induced by the Poisson structure: on each stratum $T^*(M_{H}/G)$
we have the cotangent Lie algebroid structure associated with the Poisson structure on $M_{(H)}/G$. We can also construct, via the path method, a stratified symplectic groupoid
$\Sigma(X)=\Sigma(M/G)$. In general, just like in the smooth case, this will not coincide
with the reduced symplectic groupoid $\Sigma(M)/\!/G$. The same conditions as in the smooth case apply here.

%%%%%%%%%%%%%%%%%%%%%%%%%%%%%%%%%%%%
%%%%%%%%%%%%%%%%%%%%%%%%%%%%%%%%%%%%
\section{Poisson orbispaces}       %
\label{sec:orbispaces}             %
%%%%%%%%%%%%%%%%%%%%%%%%%%%%%%%%%%%%
%%%%%%%%%%%%%%%%%%%%%%%%%%%%%%%%%%%%
%%%%%%%%%%%%%%%%%%%%%%%%%%%%%%%%%%%%

The results of the prior sections can be understood more 
geometrically in the context of what may be called \emph{Poisson orbispaces}. 
These structures  give an alternative method to deal 
with Poisson quotients avoiding singular non-smooth spaces.

%%%%%%%%%%%%%%%%%%%%%%%%%%%%%%%%%%%%%%%%%%%%%%%%
\subsection{The notion of a Poisson orbispace} %
\label{subsec:orbispace:defn}                  %
%%%%%%%%%%%%%%%%%%%%%%%%%%%%%%%%%%%%%%%%%%%%%%%%

Let $(M,\{\cdot,\cdot\},G)$ be a proper Poisson $G$-space, i.e., a Poisson
manifold $(M,\{\cdot,\cdot\})$ together with a proper Poisson action
$G\times M\to M$. The orbit space $M/G$ will be our model of a Poisson orbispace. 
However, it is convenient to adopt a slightly more general point of view. 

\begin{defn}
\label{defn:Poisson:orbifold}
Two proper Poisson $G$-spaces $G_1\times M_1\to M_1$ and $G_2\times M_2\to M_2$ 
are called \textbf{Morita equivalent} if there exists a two legs diagram
\[
\xymatrix{
&Q\ar[dl]_{p_1}\ar[dr]^{p_2}&\\
M_1& &\overline{M_2}}
\]
where:
\begin{enumerate}[(i)]
\item $Q$ is a Poisson manifold admitting commuting left and right
Poisson actions of $G_1$ and $G_2$;
\item The projections $p_i$ are $G_i$-equivariant Poisson surjective 
submersions, where $\overline{M_2}$ is the Poisson manifold whose bracket is minus the the given one on $ M_2$;
\item The projection $p_1:Q\to M_1$ (resp. $p_2:Q\to M_2$) is a
principal $G_2$-bundle (resp. $G_1$-bundle).
\end{enumerate}
A Morita equivalence class of proper Poisson $G$-spaces is called
a \textbf{Poisson orbispace}.
\end{defn}

It is easy to show, just like for ordinary Morita equivalence, that this
is an equivalence relation.  Also, in this definition, it is important to 
require $Q$ to be a Poisson manifold, and \emph{not} a symplectic manifold
(if $Q$ was symplectic, this binary relation would not even be reflexive). 
The following examples should convince the reader that this is a natural equivalence 
relation.

\begin{example}
Let $G\times M\to M$ be a free and proper Poisson action. Then $M$ is Morita
equivalent to the quotient $M/G$, with the trivial action of the trivial group.
In fact, we have the diagram
\[
\xymatrix{
&M\ar[dl]_{p_1}\ar[dr]^{p_2}&\\
M& &\overline{M/G}}
\]
where $p_1$ is the identity map and $p_2$ is the canonical projection. This diagram
fulfills all the conditions of Definition \ref{defn:Poisson:orbifold}. Conversely,
every proper Poisson $G$-space which is Morita equivalent to one having trivial 
group $G=\{e\}$, arises from a free and proper action. 
\end{example}

\begin{example}
Given a proper \'etale groupoid $\G\tto M$, a Poisson structure on $M$ is 
\emph{$\G$-invariant} if for every (local) bisection $b:M\to \G$ the composition 
$\t\circ b:M\to M$ is a (local) Poisson diffeomorphism. A \emph{Poisson orbifold} can be 
defined as a Morita equivalence class of proper \'etale groupoids with invariant Poisson
structures, as in Definition \ref{defn:Poisson:orbifold}.

Let $G$ be a finite group and $G\times M\to M$ a proper Poisson action. The quotient $M/G$
is an example of a Poisson orbifold. Conversely, every Poisson orbifold that can be written
as a global quotient is of this form but, in general, there will be different Morita equivalent presentations of this  kind. An important class is furnished by symplectic orbifolds, where the Poisson
structure is non-degenerate.
\end{example}

\begin{rem}
These examples suggest that Definition \ref{defn:Poisson:orbifold} of a Poisson orbispace is somewhat restrictive: instead of proper Poisson $G$-spaces, one should consider proper 
Lie groupoids $\G_i\to M_i$. However, in the non-\'etale case, 
one has the problem of defining $\G_i$-invariant Poisson structures on the spaces of units. If this can be solved in a reasonable way, all our results should be generalizable to such proper $\G$-spaces, using the the slice theorem of Zung \cite{Zung}.
\end{rem}

A \textbf{morphism} $\Psi:M\to N$ from a proper Poisson $G$-space $M$ to a proper Poisson 
$H$-space is a smooth map such that:
\begin{enumerate}
\item[(i)] $\Psi$ is $\psi$-equivariant relative to a homomorphism $\psi:G\to H$, that is, 
\[ \Psi(g\cdot m)=\psi(g)\cdot f(m) \quad \text{for all} \quad f \in  G, \; m \in  M;\]
\item[(ii)] $\Psi$ is a Poisson map, that is, 
\[ \{f_1\circ\Psi,f_2\circ\Psi\}_M=\{f_1,f_2\}_N\circ\Psi \quad \text{for all} \quad f_1, f _2 \in C ^{\infty}(N).\]
\end{enumerate}
A \textbf{morphism between two Poisson orbispaces} is a morphism $\Psi:M\to N$ between
any two proper Poisson $G$-spaces representing the Poisson orbispaces.

\begin{ex}
Let $M$ be a proper Poisson $G$-space and let $N\hookrightarrow M$ be any Poisson
submanifold of $M$ which is $G$-invariant. Then $N$ is a Poisson $G$-space and the inclusion
$i:N\to M$ defines a morphism of the corresponding Poisson orbispaces. Conversely, one can
define a subspace of a Poisson orbispace to be a morphism of Poisson orbispaces which
can be represented in this way.
\end{ex}

Many other constructions in Poisson geometry can be carried out in the category of Poisson 
orbispaces.

%%%%%%%%%%%%%%%%%%%%%%%%%%%%%%%%%%%%%%%%%%%%%%%%%%%%%%%%%%%%%%%%
\subsection{Poisson stratification of a Poisson orbispace}     %
\label{subsec:orbispace:quotients}                  				   %
%%%%%%%%%%%%%%%%%%%%%%%%%%%%%%%%%%%%%%%%%%%%%%%%%%%%%%%%%%%%%%%%

The following result shows that we can associate a Poisson stratification to
every orbispace.

\begin{thm}
Let $(M_1,\{\cdot,\cdot\}_1,G_1)$ and $(M_2,\{\cdot,\cdot\}_2,G_2)$ be
Morita equivalent proper Poisson spaces. Then their orbit spaces
$M_1/G_1$ and $M_2/G_2$ have isomorphic Poisson stratifications: there
exists a stratified isomorphism $\phi:M_1/G_1\to M_2/G_2$ which is also 
a Poisson map.
\end{thm}

\begin{proof}
Let $\xymatrix{M_1&Q\ar[l]_{p_1}\ar[r]^{p_2}& M_2}$ be the equivalence bimodule. 
The isomorphism $\phi$ is built in the usual way: given $[m]\in M_1/G_1$ define
$\phi([m])\in M_2/G_2$ by $\phi([m])=[p_2(q)]$, where $q\in Q$ is
any element such that $p_1(q)=m$. One checks easily that this map is well-defined 
(does not depend on choices) and gives a homeomorphism $\phi:M_1/G_1\to M_2/G_2$
which is also a Poisson map and makes the following diagram commute:
\[ 
\xymatrix{
       & Q \ar[dl]_{p_1}\ar[dr]^{p_2} \\
M_1\ar[d]_{q_1} &  & M_2\ar[d]^{q_2} \\
M_1/G_1\ar[rr]_{\phi}& & M_2/G_2 \ar@{-->}[r]_f&\Rr
}
\]
This map also preserves the stratification by
orbit types and restricts to a homeomorphism between strata. It remains to check
that it is a smooth map and that the restriction to each stratum is a diffeomorphism.

In order to check that $\phi:M_1/G_1\to M_2/G_2$ is smooth, let $f:M_2/G_2\to\Rr$ be a smooth map. 
We need to show that $f\circ\phi:M_1/G_1\to\Rr$ is smooth, i.e., that $f\circ\phi\circ q_1:M_1\to\Rr$
is smooth. But this follows by observing that, since the diagram above commutes, this map is the 
quotient by $p_1$ of the smooth map $f\circ q_2\circ p_2:Q\to\Rr$.

Let $(M_1)_{(H_1)}/G_1$ be a stratum of $M_1$ which is mapped to the stratum
$(M_2)_{(H_2)}/G_2$ of $M_2$. Then the equivalence bimodule restricts to an
equivalence bimodule
\[
\xymatrix{
&Q_{(H_1\times H_2)}\ar[dl]_{p_1}\ar[dr]^{p_2}&\\
(M_1)_{(H_1)}\ar[d]& & (M_2)_{(H_2)}\ar[d]\\
(M_1)_{(H_1)}/G_1\ar@{<->}[rr]^{\phi}& &(M_2)_{(H_2)}/G_2
}
\]
where the vertical arrows are the canonical projections to the quotients.
In this diagram all spaces are smooth manifolds and all maps (except $\phi$)
are smooth submersions. Since $\phi$ is a homeomorphism, it follows that
$\phi$ is a diffeomorphism.
\end{proof}

%%%%%%%%%%%%%%%%%%%%%%%%%%%%%%%%%%%%%%%%%%%%%%%%%%%%%%%%%%%%%%%%
\subsection{Lie theory of a Poisson orbispace}        %
\label{subsec:orbispace:Lie}                  				     %
%%%%%%%%%%%%%%%%%%%%%%%%%%%%%%%%%%%%%%%%%%%%%%%%%%%%%%%%%%%%%%%%

Every Poisson orbispace has an associated Poisson stratification, so it follows
from Section \ref{sec:symp:reduction}, that it has a stratified Lie algebroid, which
can be integrated to a stratified symplectic groupoid (under appropriate integrability
assumptions). If $(M,G)$ is any representative of the Poisson orbispace, then its stratified
symplectic groupoid is represented by $\Sigma(M/G)$ and can be obtained by integration of the stratified Lie algebroid $T^*(M/G)$ (considering equivalence
classes of cotangent paths). $\Sigma(M/G)$ is the unique source 1-connected stratified 
symplectic groupoid integrating the Poisson stratification.

Note that for any representative $(M,G)$ of a Poisson orbispace, we also 
have the stratified symplectic groupoid $\Sigma(M)/\!/G$. As we saw in Section \ref{subsec:free:reduction}, even in the free case, these two groupoids do not coincide, in general. The symplectic groupoid $\Sigma(M)/\!/G$ depends on the particular 
representation of the Poisson orbispace and so by varying the representation 
we obtain different symplectic groupoids integrating the Poisson orbispace. 

It is also natural to consider the geometric quotients $J^{-1}(0)/G$ associated 
with a representative $(M,G)$ of the orbispace (e.g., if one does not want to 
consider stratified objects). Before we look at $J:\Sigma(M)\to\gg^*$, let 
us consider first its infinitesimal version $j:T^*M\to\gg^*$ (the momentum map
for the cotangent lifted action). Since the action is Poisson, this is a Lie 
algebroid morphism, so we can think of $j^{-1}(0)\subset T^*M$ as a Lie subalgebroid.
The problem of course is that, in general, $j^{-1}(0)$ is not a smooth subbundle. 
We can still think of $j^{-1}(0)/G$ as a Lie algebroid over $M/G$ as we now explain.

Let us introduce the spaces of \textbf{basic 1-forms}, \textbf{projectable vector fields}, and  \textbf{basic vector fields} by
\begin{align*}
\Omega^1_{\rm bas}(M)&:=\{\al\in\Omega^1(M)~|~j(\al)=0\text{ and }\al\text{ is $G$-invariant}\};\\
\X_{\rm proj}(M)&:=\{X\in\X(M) \mid X(f)\in C^\infty(M)^G\text{ whenever }f\in C^\infty(M)^G\};\\
\X_{\rm bas}(M)&:=\X_{\text{proj}}(M)/\rho(C^\infty(M,\gg))), \quad \text{where} \quad \rho:\gg\to\X(M) \quad \text{is the action map}.
\end{align*}
Let us also call the $G$-invariant functions \textbf{basic functions} and denote them by $C^\infty_{\rm bas}(M)$. In the free case, we can identify basic functions, basic 1-forms, and basic vector fields with functions, 1-forms, and vector fields on $M/G$. Let us show now that in the non-free case we can still take the basic functions, the basic 1-forms, and the basic vector fields as the ingredients for the ``Lie algebroid'' of $M/G$.

The usual usual Lie bracket of vector fields gives a Lie bracket on 
$\X_{\rm bas}(M)$. On the other hand, the Lie algebroid bracket on $T^*M$ gives rise to a Lie bracket on $\Omega^1_{\rm bas}(M)$ (note that we can identify the basic 1-forms with the smooth 
sections of $j^{-1}(0)/G\to M/G$). Our next proposition shows that the basic 1-forms are the sections of 
an algebraic analog of a Lie algebroid over $M/G$, sometimes called a \emph{Lie pseudo-algebra} (see \cite{Mack}).

\begin{prop}
Let $(M,\{~,~\})$ be a proper Poisson $G$-space. The basic 1-forms, basic vector fields, and basic functions
form a Lie pseudo-algebra with Lie bracket $[~,~]$ on 
$\Omega^1_{\rm bas}(M)$ and 
anchor $\sharp:\Omega^1_{\rm bas}(M)\to \X_{\rm bas}(M)$ induced from the Lie bracket
and anchor on $TM$.
\end{prop}

\begin{proof}
It is immediate from the definitions that the Lie algebroid bracket restricts to a Lie bracket
on $\Omega^1_{\rm bas}(M)$. It is also easy to check that the anchor maps basic 1-forms to
projectable vector fields, so gives a $\Rr$-linear map $\sharp:\Omega^1_{\rm bas}(M)\to\X_{\rm bas}(M)$. We still need to check the 
Leibniz identity
\[ [\al,f\be]=f[\al,\be]+\sharp\al(f)\be,\quad \text{for all} \quad f\in C^\infty_{\rm bas}(M),\;  \al,\be\in\Omega^1_{\rm bas}(M).\]
This follows from the Leibniz identity for $T^*M$ and the fact that for any basic vector
field $[X]\in\X_{\rm bas}(M)$ with representative $X\in\X_{\rm proj}(M)$, we can define
on basic functions $[X](f):=X(f)$ which does not depend on the choice of the representatives.

\end{proof}

Finally, the pseudo-Lie algebra of a Poisson $G$-space is, in fact, 
an object naturally associated with its Poisson orbispace. This is the content of our 
next result whose proof is elementary and so can be left to the reader:
\vskip 10 pt

\begin{thm}
Let $(M_1,\{\cdot,\cdot\}_1,G_1)$ and $(M_2,\{\cdot,\cdot\}_2,G_2)$ be
Morita equivalent proper Poisson spaces. Then: 
\begin{enumerate}[(i)]
\item There are natural identifications 
\[
C^\infty_{\rm bas}(M_1)\simeq C^\infty_{\rm bas}(M_2), \quad 
\X_{\rm bas}(M_1)=\X_{\rm bas}(M_2).
\]
\item There is a natural homeomorphism of stratified spaces
\[ j_1^{-1}(0)/G\simeq j_2^{-1}(0)/G\]
which induces an isomorphism $\Omega^1_{\rm bas}(M_1)\simeq\Omega^1_{\rm bas}(M_2)$.
\end{enumerate}
These identifications give an isomorphism of the associated Lie pseudo-algebras.
\end{thm}

What about the global objects? Given Morita equivalent proper Poisson spaces $(M_1,\{\cdot,\cdot\}_1,G_1)$ and $(M_2,\{\cdot,\cdot\}_2,G_2)$, the 
corresponding group\-oid-like objects $J_1^{-1}(0)/G_1$ and $J_2^{-1}(0)/G_2$ \emph{are not} isomorphic, in general. For example, if the action is proper and free, we have that $(M,G)$
is always Morita equivalent to $(M/G,\{e\})$. However, the groupoid $\Sigma(M)/\!/G=J^{-1}(0)/G$ need not be source 1-connected, while $\Sigma(M/G)$ always is. It seems that to treat the global objects one must be willing to identify groupoids up to covers.

%%%%%%%%%%%%%%%%%%%%%%%%%%%%%%%%%%%%%%%%%%%%%%%%%%%%%%%%%%%%%%%%%%%%%%

\appendix                     

%%%%%%%%%%%%%%%%%%%%%%%%%%%%%%%%%%%%
%%%%%%%%%%%%%%%%%%%%%%%%%%%%%%%%%%%%
%%%%%%%%%%%%%%%%%%%%%%%%%%%%%%%%%%%%
\section{Poisson-Dirac and Lie-Dirac submanifolds}%
\label{appendix:submanifolds}      %
%%%%%%%%%%%%%%%%%%%%%%%%%%%%%%%%%%%%
%%%%%%%%%%%%%%%%%%%%%%%%%%%%%%%%%%%%
%%%%%%%%%%%%%%%%%%%%%%%%%%%%%%%%%%%%

\noindent {\bf Poisson-Dirac submanifolds.} 
Let $M$ be a Poisson manifold. Recall that a Poisson submanifold
$N\subset M$ is a submanifold which has a Poisson bracket and for
which the inclusion $i:N\hookrightarrow M$ is a Poisson map, that is, 
\[ 
\{f\circ i,g\circ i\}_M=\{f,g\}_N\circ i,\qquad \forall f,g\in
C^\infty(N).
\]
For example, a symplectic submanifold of a symplectic manifold is a 
Poisson submanifold if and only if it is an open subset. In fact, all 
Poisson submanifolds are unions of open subsets of symplectic leaves of $M$, which 
makes this notion too restrictive. 
For many purposes (in particular those for this paper) we need a more general notion of sub-object
in the Poisson category. These are the \emph{Poisson-Dirac}
submanifolds and in the following paragraphs we collect a few useful
facts about these submanifolds that are used in this paper. For the
proofs see \cite{CrFe1}.

\begin{defn}
\label{defn:DiracSub}
Let $M$ be a Poisson manifold with associated symplectic foliation $\F$. A submanifold $N\subset M$ is called a
\textbf{Poisson-Dirac submanifold} if $N$ is a Poisson manifold such that 
\begin{enumerate}
\item[(i)] the symplectic foliation of $N$ is $N\cap{\mathcal F}=\{L\cap
    N \mid L\in{\mathcal F}\}$, and 
\item[(ii)] for every leaf $L\in{\mathcal F}$, $L\cap N$ is a symplectic
  submanifold of $L$.
\end{enumerate}
\end{defn}

Note that if $(M,\{\cdot,\cdot\})$ is a Poisson manifold, then the
symplectic foliation with the induced symplectic forms on the leaves,
gives a smooth (singular) foliation with a smooth family of symplectic
forms. Conversely, given a manifold $M$ with a foliation ${\mathcal F}$,
by a \textbf{smooth family of symplectic forms} on the
leaves we mean a family of symplectic forms
$\set{\omega_L\in\Omega^2(L) \mid L\in\F}$ such that for every smooth
function $f\in C^\infty(M)$ the Hamiltonian vector field $X_f$ defined
by
\[ i_{X_f}\omega_L=d(f|_L),\quad \forall L\in\F\]
is a smooth vector field on $M$. Given such a smooth family, one sets $\{f,g\} := X_f(g)$.
This defines a Poisson bracket on $M$ whose 
associated symplectic foliation is precisely ${\mathcal F}$. Hence, a 
Poisson structure can be defined by specifying its symplectic foliation.  
It follows that a submanifold $N$ of a Poisson
manifold $M$ has at most one Poisson structure satisfying conditions
(i) and (ii) above and that this Poisson structure on $N$ is completely
determined by the Poisson structure of $M$ (see \cite[Section 8]{CrFe1}).

\begin{ex}
If $M$ is symplectic, then there is only one symplectic
leaf. The Poisson-Dirac submanifolds are precisely the symplectic
submanifolds of $M$.
\end{ex}

Therefore, we see that the notion of a Poisson-Dirac submanifold
generalizes to the Poisson category the notion of a symplectic submanifold.

\begin{ex}
Let $L$ be a symplectic leaf of a Poisson manifold and $N\subset M$ a
submanifold which is transverse complement to $L$ at some $x_0$, that is,
\[ T_{x_0}M=T_{x_0}L\oplus T_{x_0}N.\]
Then one can check that conditions (i) and (ii) in Definition
\ref{defn:DiracSub} are satisfied in some open subset in $N$
containing $x_0$. In other words, if $N$ is a small enough transverse
complement, then it is a Poisson-Dirac submanifold. Sometimes one calls the
Poisson structure on $N$ the transverse Poisson structure to $L$ at
$x_0$ (up to Poisson diffeomorphisms, this structure does not depend
on the choice of the transverse complement $N$).
\end{ex}

A useful property of Poisson-Dirac submanifolds is
given in the following proposition. Its proof is a simple consequence of
the definitions.

\begin{prop}
\label{prop:Poisson:Dirac:maps}
Let $N\subset M$ be a Poisson-Dirac submanifold of  $(M,\{\cdot,\cdot\})$. If 
$\phi:M\to M$ is a Poisson diffeomorphism that leaves $N$ invariant 
then its restriction $\phi|_N:N\to N$ is a Poisson diffeomorphism of $N$. 
\end{prop}

The two conditions in Definition \ref{defn:DiracSub} are not
practical to use. We give now some alternative criteria to
determine if a given submanifold is a Poisson-Dirac submanifold.

Observe that condition (ii) in the definition means that the
symplectic form on a leaf $L\cap N$ is the pull-back $i^*\omega_L$,
where $i:N\cap L\hookrightarrow L$ is the inclusion into a leaf and
$\omega_L\in\Omega^2(L)$ is the symplectic form. Let
$\sharp:T^*M\to TM$ be the bundle map determined by the Poisson structure.  We conclude that we must have
\begin{equation}
\label{eq:Dirac}
TN\cap \sharp(TN^0)=\{0\},
\end{equation}
since the left-hand side is the kernel of the pull-back $i^*\omega_L$ (for a vector subbundle $E\subset F$, we denote by $E^0\subset F^*$ its annihilator subbundle). If this condition holds, then at each point $x\in N$ we
obtain a bivector $\pi_N(x)\in\wedge^2 T_xN$ and one has the following result (see \cite{CrFe1}). 

\begin{prop}
\label{prop:Dirac}
Let $N$ be a submanifold of a Poisson manifold $M$ such that
\begin{itemize}
\item[{\rm (a)}] equation \eqref{eq:Dirac} holds, and
\item[{\rm (b)}] the induced tensor $\pi_N$ is smooth.
\end{itemize}
Then $\pi_N$ is a Poisson tensor on $N $ and $N$ is a Poisson-Dirac submanifold of $M $.
\end{prop}

Notice that, by the remarks above, the converse of the proposition also
holds.

\begin{remark}
Equation (\ref{eq:Dirac}) can be interpreted in
terms of the Dirac theory of constraints. This is the reason for the
use of the term ``Poisson-Dirac submanifold''. We refer the reader to
\cite{CrFe1} for more explanations.
\end{remark}

\medskip

\noindent {\bf Lie-Dirac submanifolds.}
From Proposition \ref{prop:Dirac} we deduce the following sufficient
condition for a submanifold to be a Poisson-Dirac submanifold.

\begin{cor}
\label{cor:LieDirac}
Let $M$ be a Poisson manifold and $N\subset M$ a submanifold. Assume
that there exists a subbundle $E\subset T_NM$ such that: 
\[ T_N M=TN\oplus E\] 
and $\sharp(E^0)\subset TN$. Then $N$ is a Poisson-Dirac submanifold.
\end{cor}

There are Poisson-Dirac submanifolds which do not satisfy the
conditions of this corollary. Also, the bundle $E$ may not be unique. 
For a detailed discussion and examples we refer to \cite{CrFe1}. 

Under the assumptions of the corollary, the Poisson bracket on the
Poisson-Dirac submanifold $N\subset M$ is quite simple to describe.
Given two smooth functions $f,g\in C^\infty(N)$, to obtain their
Poisson bracket we pick extensions $\widetilde{f},\widetilde{g}\in
C^\infty(M)$ such that $\d_x\widetilde{f}, \d_x\widetilde{g}\in
E_x^0$. Then the Poisson bracket on $N$ is given by
\begin{equation}
\label{eq:PoissonBracket}
\left\{f,g\right\}_N={\{\widetilde{f},\widetilde{g}\}|}_N.
\end{equation}
It is not hard to check that this formula does not
depend on the choice of extensions.

\begin{defn}
\label{defn:LieDirac}
Let $M$ be a Poisson manifold. A submanifold $N\subset M$ is
called a \textbf{Lie-Dirac submanifold} if there exists a subbundle $E\subset T_NM$ such that $E^0$ is a \emph{Lie subalgebroid} of $T^*M$ (equivalently, $E$ is a coisotropic submanifold of the tangent Poisson manifold $TM$).
\end{defn}

Since $E$ in the definition satisfies the assumptions of
Corollary \ref{cor:LieDirac}, Lie-Dirac submanifolds form a special class of
Poisson-Dirac submanifolds with very special geometric
properties, first studied by Xu \cite{Xu}. For example,
a Poisson-Dirac submanifold of an integrable Poisson manifold may not
be integrable, but for Lie-Dirac submanifolds we have (see \cite{Xu,CrFe1}) the following result.

\begin{prop}
\label{prop:Lie:Dirac}
Let $N\subset M$ be a Lie-Dirac submanifold of $(M,\{\cdot,\cdot\})$. 
If $M$ is an integrable Poisson manifold, then $N$ is also integrable
and $\Sigma(N)$ is a symplectic subgroupoid of $\Sigma(M)$.
\end{prop}

In fact, $\Sigma(N)$ is the Lie subgroupoid which corresponds to the
Lie subalgebroid $E^0\subset T^*M$.

% -----------------------------------------------------------------------

% -----------------------------------------------------------------------
\end{document}